\documentclass[epic,eepic,11pt]{amsart}

\usepackage{amsfonts,mathrsfs}
\usepackage[mathscr]{eucal}
\usepackage{mathtools}
\usepackage{amssymb}
\usepackage{amscd}
\usepackage{fancyhdr}
\usepackage[all,cmtip]{xy}

 \usepackage[all]{xy}
\usepackage{tikz-cd}

\usepackage{euler,eucal}

\DeclareMathAlphabet{\mathbf}{T1}{ppl}{bx}{n}
\DeclareMathAlphabet{\mathrm}{T1}{ppl}{m}{n}

\numberwithin{equation}{section}

\newcommand\note[1]%
{$^\dagger$\marginpar{\footnotesize{$^\dagger${#1}}}}

\def\({\left(}
\def\){\right)}
\def\<{\left<}
\def\>{\right>}

\newtheorem{theorem}{Theorem}[section]
\newtheorem{proposition}[theorem]{Proposition}
\newtheorem{lemma}[theorem]{Lemma}
\newtheorem{definition}[theorem]{Definition}

\newtheorem{corollary}[theorem]{Corollary}

\theoremstyle{definition}

\newtheorem{remark}[theorem]{Remark}

\newcommand\R{\mathbb{R}}

\newcommand     {\comment}[1]   {}
\newcommand{\mute}[2] {}
\newcommand     {\printname}[1] {}

\newcommand\func[1]{\operatorname{\mathrm{#1}}}
\newcommand\funclim[1]{\operatorname*{\mathrm{#1}}}

\renewcommand\dim{\func{dim}}

\renewcommand\ker{\func{ker}}
\renewcommand\lim{\funclim{lim}}
\renewcommand\log{\func{log}}

\newcommand\sur{\mathrel{\to\kern-1.8ex\to}}
\newcommand\iso{\mathrel{\hookrightarrow\kern-1.8ex\to}}

\newcommand\longhookrightarrow{\lhook\joinrel\longrightarrow}

\newcommand\longsur{\mathrel{\longrightarrow\kern-1.8ex\to}}
\newcommand\longiso{\mathrel{\longhookrightarrow\kern-1.8ex\to}}

\renewcommand\supset{\supseteq}
\renewcommand\emptyset{\varnothing}

\begin{document}

\title{Equivariant Morse theory for Lie algebra actions on Riemannian foliations}

\author{Yi Lin, Zuoqin Wang}

\thanks{Z.W. is partially supported by  National Key R and D Program of China 2020YFA0713100, and by NSFC no. 12171446.}

\address{Yi Lin \\ Department of Mathematical Sciences \\  Georgia Southern University \\Statesboro, GA, 30460 USA}
\email{yilin@georgiasouthern.edu}

\address{Zuoqin Wang \\ School of Mathematical Sciences \\  University of Science and Technology of China \\Hefei, Anhui, 230026 P.R.China}

\email{wangzuoq@ustc.edu.cn}

\date{\today}


\begin{abstract}
We consider a transverse isometric action of a finite-dimensional Lie algebra
$\mathfrak g$ on a Riemannian foliation. In this setting, we study
equivariant Morse--Bott theory on the leaf space of the foliation. Among
other results, we establish a foliated Morse--Bott lemma for
$\mathfrak g$-invariant basic Morse--Bott functions and a foliated analogue
of the usual handle presentation theorem. In the non-equivariant case, we use
these results to give a new proof of the Morse inequalities for Riemannian
foliations. In the equivariant case, we apply them to Hamiltonian actions of
abelian Lie algebras on presymplectic manifolds whose underlying foliations
are Riemannian, and we extend the Kirwan surjectivity and injectivity
theorems from equivariant symplectic geometry to this setting. As a
consequence, Kirwan surjectivity and injectivity hold for Hamiltonian torus
actions on symplectic orbifolds.
\end{abstract}

\maketitle
\tableofcontents

\section{Introduction}

Assume that there is a Hamiltonian action of a compact Lie group on a symplectic manifold $(M, \omega)$ with a moment map
$\mu: M\rightarrow \mathfrak{g}^*$ taking values in the dual of the Lie algebra of $G$, and that $0\in \mathfrak{g}^*$ is a regular value of $\mu$. Using the equivariant Morse theory, Kirwan \cite{Kir84} proved that the map $\kappa : H_G(M) \rightarrow H_G(\mu^{-1}(0))$ induced by inclusion is a surjection, and that the pullback homomorphism $i^*: H_G(M)\rightarrow H_G(M^G)$ is injective when $G$ is a compact and connected torus. These two results, known respectively as the Kirwan surjectivity and
Kirwan injectivity theorems, have played an important role in the development of equivariant symplectic geometry. More recently, Lin and Sjamaar  \cite{LS17} studied Hamiltonian actions on presymplectic manifolds. They discovered that when the action is clean, components of a moment map are Morse--Bott functions, and extended the Atiyah-Guillemin-Sternberg convexity theorem to this setting. It is natural to ask whether the Kirwan surjectivity and injectivity theorems could also be generalized to Hamiltonian actions on presymplectic manifolds. In \cite{LY19} and \cite{LY23}, using symplectic Hodge theoretic techniques, Lin and Yang generalized the Kirwan injectivity theorem to the case of a Hamiltonian torus action on a presymplectic manifold that satisfies the transverse Hard Lefschetz property.

The notion of basic cohomology was first introduced by Reinhart \cite{R59} as a cohomology theory for the leaf space of a foliation. It turns out to be a very useful cohomological tool in the study of Riemannian foliations. Killing foliation is an important class of Riemannian foliations which is known to possess a type of "internal'' symmetry given by the transverse action of their structural Lie algebras. In order to develop the localization techniques for this important type of symmetries, Goertsches and T\"{o}ben \cite{GT18}  first came up with the notion of equivariant basic cohomology. In particular, Reeb flow of a $K$-contact manifold is isometric with respect to the given contact metric and so the underlying foliation is Killing. For the action of a torus on a $K$-contact manifold,
 Casselmann \cite{Ca17} applied the machinery of equivariant Morse theory developed in \cite{W69} to study its equivariant basic cohomology, and proved that the Kirwan surjectivity and injectivity continue to hold in this framework. 


The theory of complete pseudogroups of local isometries was first introduced by A Haefliger in 1980's as an important approach to the transverse geometry and topology of Riemannian foliations. From this viewpoint,  L$\acute{o}$pez \cite{AL93} studied the Morse theory for the orbit space of complete pseudogroup of local isometries. When the orbit space is compact,  using Witten's analytic methods L$\acute{o}$pez established Morse inequalities for the invariant cohomology in this setting. Among other things, his result implies the Morse inequalities for the basic cohomology of Riemannian foliations with compact leaf spaces. 


In this paper, we initiate the study of equivariant Morse theory on a general Riemannian foliation. Suppose that there is an isometric transverse action of a Lie algebra $\mathfrak{g}$ on a Riemannian foliation $(M, \mathcal{F})$. We prove the following three Morse-theoretic results.
Assume that $f:M\rightarrow \mathbb R$ is a proper
$\mathfrak g$-invariant basic Morse--Bott function on $(M,\mathcal F)$
which attains a minimum. We prove the following three Morse-theoretic
results.

\begin{lemma}(\textbf{Foliated Morse--Bott Lemma})
Let $X$ be a connected, compact and nondegenerate critical submanifold of
$f$. For any $p\in X$, define $N_p^+X$ and $N_p^-X$ to be the positive and
negative subspaces of the quadratic form
$T_p^2f:N_pX\times N_pX\rightarrow \mathbb R$. Then both
$N^+X:=\cup_{p\in X}N_p^+X$ and
$N^-X:=\cup_{p\in X}N_p^-X$ are subbundles of $NX$ such that
\begin{equation}
NX=N^+X\oplus N^-X.
\end{equation}

Moreover, if $f(X)=0$, then there exist a
$(\mathfrak g\ltimes\mathcal F)$-equivariant tubular neighborhood
$\phi:NX\rightarrow M$, a positive constant $r>0$, and a
$(\mathfrak g\ltimes\mathcal F)$-equivariant diffeomorphism $\theta$ from
$NX(r)$ to an open neighborhood $U$ of the zero section in $NX$, which is
fiberwise and origin preserving, such that
\begin{equation}
f((\phi\circ\theta)(x,y))
=
\|x\|^2-\|y\|^2,
\qquad
\forall\, (x,y)\in N_pX,\quad \text{where } p\in X.
\end{equation}
\end{lemma}

\begin{theorem}
Suppose that $a<b$ are two real numbers in $f(M)$, and that no critical
values of $f$ lie in the open interval $(a,b)$. Then there is a
foliation-preserving $\mathfrak g$-equivariant diffeomorphism from $M^b$
onto $M^a$.
\end{theorem}

\begin{theorem}
Suppose that $0$ is a critical value of $f$, that $\epsilon>0$ is so small
that $0$ is the only critical value of $f$ in $(-\epsilon,\epsilon)$, and
that $X_1,\ldots, X_s$ are all the connected components of the critical
submanifold of $f$ contained in $f^{-1}(0)$. Then there is a
$\mathfrak g$-equivariant foliated diffeomorphism from $M^{\epsilon}$ to a
$(\mathfrak g\ltimes\mathcal F)$-invariant submanifold
$\tilde M\subset M$ with boundary, obtained from $M^{-\epsilon}$ by attaching
$s$ basic handles of type $(N^+X_i,N^-X_i)$, $1\leq i\leq s$.
\end{theorem}

These are Lemma \ref{f-Morse--Bott}, Theorem \ref{regular-values} and Theorem \ref{crossing-critical} in this article, respectively.  There are certain technical subtleties in extending Morse theory to Riemannian foliations.  For instance, let $f$ be a basic Morse function on $M$ and let $M^a:=\{p\in M\,\vert\, f(p)\leq a\}$. Note that since the basic cohomology is defined using a subcomplex of smooth differential forms, the usual description of the change of continuous homotopy type of $M^a$ when passing a critical level does not apply to the equivariant basic cohomology directly. Also note that when passing a critical level the attached basic handle is by definition a manifold with corner itself. To
circumvent the technical issues of working with forms on manifolds with corners, we carefully choose a version of the relative equivariant basic cohomology which allows us to give a clean proof of the desired excision property (Proposition \ref{c-supported-cohomology}).

As a main application of the equivariant Morse theoretic results established in this paper, we study the Hamiltonian transverse Lie algebra actions on presymplectic manifolds and prove the following results.

\begin{theorem} (\textbf{Kirwan Surjectivity}) Consider the transverse isometric action of an abelian Lie algebra $\mathfrak{g}$ on a Riemannian foliation $(M, \mathcal{F})$ that is also transversely symplectic. Suppose that the action of $\mathfrak{g}$ is Hamiltonian with a moment map $\mu : M \rightarrow \mathfrak{g}^*$ for which $0 \in\mathfrak{g}^*$ is a regular value. Then the inclusion map $i: Z:=\mu^{-1}(0)\xhookrightarrow{} M$ induces a surjection $\kappa: H^*_{\mathfrak{g}}(M, \mathcal{F})\rightarrow H^*_{\mathfrak{g}}(Z, \mathcal{F}\vert_Z)$ in cohomology.

\end{theorem}

\begin{theorem}(\textbf{Kirwan Injectivity}) 
Consider a transverse isometric action of an abelian Lie algebra
$\mathfrak{g}$ on a Riemannian foliation $(M,\mathcal{F})$ that is also
transversely symplectic. Assume that the action of $\mathfrak{g}$ is
Hamiltonian with moment map $\mu:M\rightarrow\mathfrak{g}^*$. Let
$X:=M^{\mathfrak{g}}$ be the fixed-leaf set, let $i:X\hookrightarrow M$ be
the inclusion map, and let
\begin{equation} \label{pullback}
i^*: H_{\mathfrak{g}}(M,\mathcal{F})\rightarrow
H_{\mathfrak{g}}(X,\mathcal{F}\vert_X)
\end{equation}
be the pullback morphism induced by $i$. Suppose that there exists a generic
element $\xi\in\mathfrak{g}$ such that, for
$\mathfrak h:=\text{span}\{\xi\}$, one has $M^{\mathfrak h}=M^{\mathfrak g}$,
and such that the component $\mu^\xi$ is proper and attains a minimum on
$M$. Then the morphism \eqref{pullback} is injective.
\end{theorem}

These are Theorem \ref{kirwan-surj-general} and Theorem \ref{Kirwan-inj-2} in this article, respectively. It is worth mentioning that to prove the Kirwan surjectivity for an abelian Lie algebra of arbitrary dimension, we follow the idea of an inductive argument used in \cite{BL10} and reduce the general case to the case when the dimension of the Lie algebra is one.  These results are new in that they provide a general framework for the Kirwan surjectivity and injectivity theorems to continue to hold on singular spaces. For instance, it was proved by Holm and Matsumura \cite{HM12} that Kirwan injectivity theorem holds for Hamiltonian torus actions on a symplectic orbifold that can be realized as the global quotient of a smooth manifold by another torus.  Since every effective orbifold can be realized as the leaf space of a Riemannian foliation, as shown in this paper our general result implies immediately that both Kirwan surjectivity and injectivity hold for Hamiltonian torus actions on a symplectic orbifold which is not necessarily the global quotient of a smooth manifold by another torus.

It is well-known that orbifolds are singular spaces locally modeled by quotients of Cartesian spaces by the smooth action of finite groups. In order to generalize the Delzant construction in toric geometry to simple non-rational polytopes, Prato \cite{Pra01} introduced the notion of quasifolds, which are singular spaces locally modeled by quotients of Cartesian spaces by the smooth affine action of countable groups. Battaglia and Zaffran later showed in \cite{BZ15} and \cite{BZ17} that every symplectic toric quasifold is equivariantly symplectomorphic to the leaf space of a Killing foliation. In a very recent joint work Lin and Miyamoyto \cite{LM23} proved that the leaf space of a Killing foliation always carries the structure of a quasifold. In view of these results, our general results will also apply to a large class of Hamiltonian torus actions on symplectic quasifolds. Indeed, the Kirwan injectivity theorem has already been used in \cite{LY23} to derive explicit combinatorial formulae for the basic Betti numbers and basic Hodge numbers of toric quasifolds.

Finally, we note that Ortiz and Valencia \cite{OV24} recently began the
study of Morse theory on Lie groupoids, a topic closely related to ours.
However, in \cite{OV24} the Lie groupoids are required to be proper in order
for the standard machinery of Morse theory to continue to work. By contrast,
our approach builds on Molino's structure theory of Riemannian foliations and
does not require the holonomy groupoid of the foliation to be proper. This is
a significant difference, since the holonomy groupoid of a Riemannian
foliation is known to be non-proper in general. 

A key ingredient in our approach is the passage to the Molino manifold.
However, this passage does not reduce the equivariant basic cohomological
results of the present paper to the Morse-theoretic results on proper Lie
groupoids obtained in \cite{OV24}. Let $\pi:P\to M$ denote the transverse
orthonormal frame bundle of $(M,\mathcal F)$, and let $\rho:P\to W$ denote
the Molino fibration, whose fibers are the leaf closures of the lifted
foliation on $P$. An isometric transverse $\mathfrak g$-action on
$(M,\mathcal F)$ can be transported through the Molino bi-bundle
$M\xleftarrow{\pi}P\xrightarrow{\rho}W$ to a genuine isometric
$\mathfrak g$-action on the Molino manifold $W$. When $W$ is compact, the
closure of the corresponding group action in $\operatorname{Iso}(W)$ is
compact, and hence gives rise to a proper action groupoid.

Nevertheless, this construction records only the geometry of leaf closures.
It does not retain the full leafwise information needed to construct
foliation-preserving diffeomorphisms on $(M,\mathcal F)$. In addition, the
basic de Rham complex of $(M,\mathcal F)$ is not naturally identified, via
the Molino construction, with the de Rham complex of $W$. If $\alpha$ is a
basic form on $(M,\mathcal F)$, then $\pi^*\alpha$ is basic for the lifted
foliation on $P$, but it need not be basic with respect to the leaf-closure
foliation on $P$ and therefore need not descend to $W$. Thus Morse-theoretic
results on the proper groupoid associated with the compact group action on
$W$ do not, by themselves, imply the equivariant basic cohomology results
established in this paper.

In the special case of orbifolds, proper Lie groupoids enter naturally, and
the Morse-theoretic results of \cite{OV24} are closely related to the
nonequivariant Morse-theoretic part of the present paper. However,
\cite{OV24} does not establish the particular equivariant cohomological
ingredients needed for the Kirwan surjectivity and injectivity theorems
proved here, such as the equivariant Thom isomorphism and the Euler-class
non-zero-divisor arguments in the equivariant Cartan-model setting. Even in
classical equivariant symplectic geometry, the Kirwan package is not a formal
consequence of Morse--Bott theory alone; it also uses equivariant
cohomological input such as Thom isomorphisms and Euler-class arguments. The
present paper develops the foliated analogues of these ingredients and
applies them to Hamiltonian torus actions on symplectic orbifolds.
 
 Our paper is organized as follows. Section \ref{review} reviews basic notions on foliations, especially Molino's structure theory of a Riemannian foliation. Section \ref{tran-action} reviews background materials on transverse Lie algebra actions on a foliated manifold. Section \ref{eq-cohomology} reviews the equivariant basic cohomology and the equivariant Thom isomorphism in this setting. 
 Section \ref{Eq-Morse} extends the equivariant Morse theory to the setting of isometric transverse Lie algebra actions on Riemannian foliations. Section \ref{application1} discusses the applications of the Morse theory to the equivariant basic cohomology on a Riemannian foliation.  Section \ref{Morse-inequalities} derives the Morse inequality for the basic cohomology of a Riemannian foliation. Section \ref{pre-sympl} generalizes the Kirwan surjectivity and Kirwan injectivity theorem to Hamiltonian actions of abelian Lie algebras on presymplectic manifolds. Section \ref{orbifold} proves the Kirwan surjectivity and Kirwan injectivity results for Hamiltonian torus actions on symplectic orbifolds.

\section{Review of Riemannian foliations}\label{review}

In this section we review the definition of a Riemannian foliation, as well
as the main ingredients of Molino's structure theory of Riemannian foliations
that will be used in this paper. We refer to \cite{Mo88} and
\cite[Sec.~3]{LS18} for a more detailed account. We begin by recalling some
basic notions in foliation theory.

\begin{definition}\label{f-chart}
Let $M$ be a manifold of dimension $p+q$. An atlas of \textbf{foliation
coordinate charts} is a maximal family of coordinate charts of the form
\[
\varphi_{\alpha}:U_{\alpha}\longrightarrow
\varphi_{\alpha}(U_{\alpha})=V_{\alpha}\times W_{\alpha}
\subset \mathbb R^p\times\mathbb R^q,
\qquad
z\longmapsto
(x_{\alpha}^1(z),\ldots,x_{\alpha}^p(z),
y_{\alpha}^1(z),\ldots,y_{\alpha}^q(z)),
\]
where $V_{\alpha}\subset \mathbb R^p$ and
$W_{\alpha}\subset \mathbb R^q$ are open subsets, such that on each overlap
$U_{\alpha}\cap U_{\beta}$ the coordinate change
$\varphi_{\beta}\circ\varphi_{\alpha}^{-1}:
\varphi_{\alpha}(U_{\alpha}\cap U_{\beta})\to
\varphi_{\beta}(U_{\alpha}\cap U_{\beta})$ satisfies
\begin{equation}\label{coordinate-change}
\dfrac{\partial y_{\beta}^i}{\partial x_{\alpha}^j}=0,
\qquad
1\leq i\leq q,\quad 1\leq j\leq p.
\end{equation}
Let $(U_{\alpha},\varphi_{\alpha},x_{\alpha}^1,\ldots,x_{\alpha}^p,
y_{\alpha}^1,\ldots,y_{\alpha}^q)$ be a foliation coordinate chart on $M$.
A \textbf{vertical plaque} in $U_{\alpha}$ is a non-empty subset of
$U_{\alpha}$ of the form
\[
y_{\alpha}^1=c_1,\ldots,y_{\alpha}^q=c_q,
\]
where $(c_1,\ldots,c_q)\in\mathbb R^q$. We denote by $(M,\mathcal F)$, or
simply by $\mathcal F$, the foliation defined by such an atlas of foliation
coordinate charts on $M$.
\end{definition}

\begin{remark}  In this paper we will have to deal with foliated manifolds with boundary. When $M$ is a bounded manifold of dimension $p+q$, replacing $\mathbb{R}^q$ by $\mathbb{R}^q_+=\{(y_1, \cdots, y_q)\in\mathbb{R}^q\,\vert\, y_q\geq 0\}$ in Definition \ref{f-chart} provides us the definition of an atlas of foliation coordinate charts on a bounded manifold.

\end{remark}

Let $\{(U_{\alpha}, x_{\alpha}^1,\cdots, x^p_{\alpha}, y_{\alpha}^1, \cdots, y_{\alpha}^q)\}$ be an atlas of foliation coordinate charts on a 
manifold $M$. Then on each $U_{\alpha}$, $\text{span}\{\dfrac{\partial}{\partial x_{\alpha}^1}, \cdots, \dfrac{\partial}{\partial x_{\alpha}^p}\}$ defines an integrable distribution; furthermore, it is easy to check that the definition is independent of the choice of a local coordinate chart, and so gives rise to an integrable distribution on $M$. A leaf of $(M, \mathcal{F})$ is a maximal integrable submanifold of this integrable distribution.

Let $\mathcal{F}$ be a foliation on a smooth manifold $M$.
Throughout this paper we denote by $\mathfrak{X}(M)$ the space of vector fields on $M$ and by
$\mathfrak{X}(\mathcal{F})\subset\mathfrak{X}(M)$ the subspace of vector fields tangent to the leaves of $\mathcal{F}$.
We say that a subset $X$ is $\mathcal{F}$-\textbf{saturated}, if whenever a leaf $L$ of $\mathcal{F}$ intersects $X$ non-trivially, we must have that $L\subset X$.
We say that a vector field $X\in\mathfrak{X}(M)$ is \emph{foliate}, if $[X,Y]\in\mathfrak{X}(\mathcal{F})$ for all $Y\in \mathfrak{X}(\mathcal{F})$.
We will denote by
$\mathfrak{R}(\mathcal{F})$ the space of foliate vector fields on $(M,\mathcal{F})$.
Clearly we have that
$\mathfrak{X}(\mathcal{F}) \subset \mathfrak{R}(\mathcal{F})$.
In this context,  a \textbf{transverse vector field} is an equivalence class in the quotient space $\mathfrak{R}(\mathcal{F})/\mathfrak{X}(\mathcal{F})$. The space of transverse vector fields, denoted by $\mathfrak{X}(M, \mathcal{F})$, forms a Lie algebra with a Lie bracket inherited from the natural one on $\mathfrak{R}(\mathcal{F})$. It is straightforward to check that $\mathfrak{X}(M, \mathcal{F})$ equals the space of $\mathfrak{X}(\mathcal{F})$-fixed sections of the normal bundle $N\mathcal{F}$ of the foliation.

The space of \textbf{basic functions}  on a foliated manifold $(M, \mathcal{F})$ is defined to be 
\[ \{ f\in C^{\infty}(M)\,\vert\, \mathcal{L}_Xf=0, \forall\, X\in \mathfrak{X}(\mathcal{F})\}.\]  More generally, 
the space of \textbf{basic forms} on $(M, \mathcal{F})$ is defined to be
$$
\Omega(M,\mathcal{F})=
\bigl\{\alpha\in\Omega(M)\,|\,
\iota(X)\alpha=\mathcal{L}(X)\alpha=0,\,
\forall
\,X\in\mathfrak{X}(\mathcal{F})\bigr\}.
$$
Since the exterior differential operator $d$ preserves basic forms, we obtain a subcomplex  $\{\Omega^{*}(M,\mathcal{F}),d\}$ of the usual de Rham complex, called the \textbf{basic de Rham complex}.
The associated cohomology $H^{*}(M,\mathcal{F})$
is called the \textbf{basic cohomology}.


\begin{definition}
A \textbf{transverse Riemannian metric} on a foliation $(M,\mathcal{F})$ is a Riemannian metric $g$ on the normal bundle $N\mathcal{F}$ of the foliation, such that $\mathcal{L}(X)g=0$,
$\forall\,X\in\mathfrak{X}(\mathcal{F})$.
We say that $\mathcal{F}$ is a \textbf{Riemannian foliation} if there exists a transverse Riemannian metric on $(M,\mathcal{F})$.
\end{definition}

\begin{definition}  We say that a Riemannian metric $g_{TM}$
on a foliated manifold $(M, \mathcal{F})$ is \textbf{bundle-like}, if it has the following property: the function $g_{TM}(X, Y)$ is basic for all foliate vector fields $X$ and $Y$ that are perpendicular to the leaves of $\mathcal{F}$.
\end{definition}
A bundle-like metric $g_{TM}$ gives rise to a transverse metric $g$ by identifying $N\mathcal{F}$ with the $g_{TM}$-orthogonal complement of
$T\mathcal{F}$ and then restricting $g_{TM}$ to $N\mathcal{F}$ . Conversely, for every transverse metric $g$ there is a bundle-like metric $g_{TM}$ which induces $g$, cf. \cite[Sec. 3.2]{Mo88}. We say a Riemannian foliation $(M, \mathcal{F})$ is \textbf{complete}, if it admits a complete bundle-like metric.

Let $(M,\mathcal{F})$ be a Riemannian foliation with a transverse metric $g$, and let $\overline{X}$ be a transverse vector field.
Define
$\mathcal{L}(\overline{X})g=\mathcal{L}(X)g$,
where $X$ is a foliate vector field that represents $\overline{X}$.
It is straightforward to check that this definition does not depend on the choice of a foliate vector representing $\overline{X}$.
A transverse vector field $\overline{X}$ is said to be \emph{transversely Killing} if $\mathcal{L}(\overline{X})g=0$.
Suppose that both $\overline{X}$ and $\overline{Y}$ are transversely Killing.
Then it follows easily from the Cartan identities that $[\overline{X},\overline{Y}]$ is also transversely Killing.
In other words, the space of transversely Killing vector fields, which we denote by $\mathfrak{X}(M, \mathcal{F}, g)$, forms a Lie subalgebra of $\mathfrak{X}(M, \mathcal{F})$.


 Let $g$ be a transverse Riemannian metric on a complete Riemannian foliation $(M, \mathcal{F})$, let $\pi: P\rightarrow M$ be the transverse orthonormal frame bundle associated to $g$, and let $K=O(q)$ be the structure Lie group of $P$.  Then the foliation $\mathcal{F}$ naturally lifts to a transversely parallelizable foliation $\mathcal{F}_P$ on $P$ that is invariant under the action of $K$; moreover, there is a locally trivial smooth fibration $\rho: P\rightarrow W$ whose fibers are leaf closures of $\mathcal{F}_P$. Here $W$ is a smooth manifold called the \textbf{Molino manifold} of the Riemannian foliation $(M, \mathcal{F})$. We say that $(M, \mathcal{F})$ is \textbf{transversely compact}, if its Molino manifold $W$ is compact. 

A transverse Riemannian metric $g_P$ on $P$ can be constructed as follows.
Since the foliation is Riemannian, the transverse Levi-Civita connection
$\theta_{LC}$ defines a basic $1$-form on $P$ with values in the Lie algebra
$k=o(q)$. Choose an invariant inner product on $k$. Let
$f_1,\ldots,f_r$ be an oriented orthonormal basis of $k^*$, and set
$\alpha_i=\theta_{LC}^*(f_i)$, $1\leq i\leq r$. Then
$\alpha_1,\ldots,\alpha_r$ are $\mathcal F_P$-basic $1$-forms, and
\begin{equation}\label{t-metric}
g_P=\pi^*g+\sum_{i=1}^r\alpha_i\otimes\alpha_i
\end{equation}
defines a transverse Riemannian metric on $P$. Similarly, if $g_{TM}$ is the
bundle-like metric on $M$ inducing $g$, then
$g_{TP}:=\pi^*g_{TM}+\sum_{i=1}^r\alpha_i\otimes\alpha_i$ is a bundle-like
metric on $P$ inducing $g_P$. It is clear that
$\pi:(P,g_{TP})\to(M,g_{TM})$ is a Riemannian submersion.

Let $X$ be a foliate vector field on $(M,\mathcal F)$. Then $X$ naturally
lifts to a foliate vector field $X_P$ on $P$. If $X$ is transversely Killing
with respect to $g$, then it follows easily from \eqref{t-metric} that
$X_P$ preserves the transverse Riemannian metric $g_P$. If $X$ is tangent to
the leaves of $\mathcal F$, then $X_P$ is tangent to the leaves of the lifted
foliation $\mathcal F_P$. These observations give rise to a Lie algebra
homomorphism
\begin{equation}\label{lifting-homo}
\pi^{\sharp}: \mathfrak X(M,\mathcal F,g)
\longrightarrow
\mathfrak X(P,\mathcal F_P,g_P).
\end{equation}

Now consider the fibration $\rho:P\to W$. Let $u$ and $v$ be vector fields on
$W$. There are unique vector fields $\tilde u$ and $\tilde v$ on $P$ which
project to $u$ and $v$, respectively, under the tangent map of $\rho$ and
which are perpendicular to the fibers of $\rho$. Since $\tilde u$ and
$\tilde v$ are projectable, they are foliate. Since they are perpendicular to
the leaf closures of $\mathcal F_P$, they are also perpendicular to the
leaves of $\mathcal F_P$. Thus $g_{TP}(\tilde u,\tilde v)$ is an
$\mathcal F_P$-basic function, and hence is also
$\overline{\mathcal F}_P$-basic. It follows that $g_{TP}$ naturally descends
to a Riemannian metric $g_W$ on $W$, such that
$\rho:(P,g_{TP})\to(W,g_W)$ is a Riemannian submersion.

In summary, we have the following Molino diagram:
\begin{equation}\label{molino-diagram}
\begin{tikzcd}
  &(P,g_{TP})  \arrow[dl, swap, "\pi"] \arrow[dr, "\rho"]&\\
  (M,g_{TM}) & &(W,g_W)
\end{tikzcd}
\end{equation}

The action of the structure Lie group $K$ on $P$ descends to an action of
$K$ on $W$ preserving the Riemannian metric $g_W$. It follows from
\eqref{molino-diagram} that $M/\overline{\mathcal F}$ is homeomorphic to
$W/K$, which is a Whitney stratified space. Moreover, a transverse Killing
vector field $X_P\in\mathfrak X(P,\mathcal F_P,g_P)$ projects under the
tangent map of $\rho$ to a Killing vector field $X_W$ on the Riemannian
manifold $(W,g_W)$. Denote by $\mathfrak X(W,g_W)$ the Lie algebra of
Killing vector fields on $(W,g_W)$. Then we have a Lie algebra homomorphism
\begin{equation}\label{transport-killing}
\rho_*\circ\pi^{\sharp}:
\mathfrak X(M,\mathcal F,g)
\longrightarrow
\mathfrak X(W,g_W).
\end{equation}

\section{Transverse Lie algebra actions} \label{tran-action}
In this section, we review necessary background materials on transverse action of a Lie algebra on a foliation. In particular, results established in \cite{LS18} on transverse isometric Lie algebra actions on Riemannian foliations play an important role in this work. 


\subsection{Transverse action Lie algebroid}

\begin{definition} (\cite{GT18})
A \emph{transverse action} of a Lie algebra $\mathfrak{g}$ on a foliated manifold $(M,\mathcal{F})$ is a Lie algebra homomorphism \begin{equation}\label{equ1.1}
a: \mathfrak{g}\rightarrow\mathfrak{X}(M, \mathcal{F}).  \end{equation}
\end{definition}

\begin{definition}(\cite{LS18})
Let $\mathfrak{g}_M$ be the trivial bundle with fiber $\mathfrak{g}$ over $M$. We define the transverse action Lie algebroid $\mathfrak{g}\ltimes \mathcal{F}$ of the transverse $\mathfrak{g}$-action to be the fibered product
\[\mathfrak{g}\!\ltimes\! \mathcal{F} \!=\! \mathfrak{g}_M \!\times_{ N\mathcal{F}}\! TM\!=\!\{(x,\zeta,v)\in M\!\times\! \mathfrak{g}\!\times\! TM\, \vert\,v\!\in\! T_x M,\, \text{and}\, \zeta_{M,x} \!=\!v\, \text{mod}\, T_x\mathcal{F}\},\]
which is a smooth subbundle of the bundle $\mathfrak{g}_M \times T M$ over $M$. Define the anchor map $t: \mathfrak{g}\ltimes \mathcal{F}\rightarrow TM$  by $t(\zeta, v)=v$.  For every smooth map $\zeta : M\rightarrow \mathfrak{g}$ and every $x\in M$ we define $\zeta_x \in N_x\mathcal{F}$ to be the value of the transverse vector field $a(\zeta) \in \mathfrak{X}(M, \mathcal{F} )$ at $x$. A smooth section of $\mathfrak{g}\ltimes \mathcal{F}$ is a pair $(\zeta,v)\in C^{\infty}(M,\mathfrak{g})\times \mathfrak{X}(M)$ satisfying $\zeta_x =v_x\, \text{mod}\,T_x\mathcal{F}$ for all $x\in M$.We define the bracket of two sections $(\zeta, v)$ and $(\eta, w)$ by
\begin{equation}\label{Lie} [(\zeta, v), (\eta, w)](x) = ([\zeta(x), \eta(x)] + L(v)(\eta)(x)-L(w)(\zeta)(x), [v, w](x))\end{equation}
\end{definition}




 Suppose that there is a transverse Lie action $a: \mathfrak{g}\rightarrow \mathfrak{X}(M,\mathcal{F})$. 
The orbit $(\mathfrak{g}\ltimes\mathcal{F})(x)$ of the transverse action Lie algebroid $\mathfrak{g}\ltimes \mathcal{F}$  is by definition a leaf through $x\in M$ of the singular foliation nduced by the image of the anchor map $t$, which we also refer to as the $\mathfrak{g}$-orbit of the leaf $\mathcal{F}(x)$. A subset $A\subset M$ is said to be $(\mathfrak{g}\ltimes \mathcal{F})$-invariant if a 
$\mathfrak{g}$-orbit of a leaf $L$ intersects $A$ non-trivially, then the entire $\mathfrak{g}$-orbit must lie inside $A$.
We say a function on $M$ is $(\mathfrak{g}\ltimes \mathcal{F})$-invariant, if for any leaf $L$,  its restriction to a $\mathfrak{g}$-orbit of $L$ is a constant. The stabilizer of $x$ for the transverse action Lie algebroid $\mathfrak{g} \ltimes \mathcal{F}$ is defined to be
\begin{equation}\label{stabilizer} \text{stab}(x, \mathfrak{g}\ltimes \mathcal{F} )= \{ \xi\in \mathfrak{g}\,\vert\, a(\xi)(x)=0\}. \end{equation}

\begin{lemma}(\cite[Lemma 2.2.7]{LS18})
Suppose that $x$ and $y$ lie in the same leaf of $\mathcal F$. Then
$a(\xi)(x)=0$ if and only if $a(\xi)(y)=0$ for every $\xi\in\mathfrak g$.
As a consequence,
\[
\text{stab}(x,\mathfrak g\ltimes\mathcal F)
=
\text{stab}(y,\mathfrak g\ltimes\mathcal F).
\]
\end{lemma}

Let $(M,\mathcal F)$ and $(M',\mathcal F')$ be two foliated manifolds, and
let $f:M\to M'$ be a foliate map. By assumption, for every
$X_p\in T_p\mathcal F$, one has $f_{*,p}(X_p)\in T_{f(p)}\mathcal F'$.
Thus $f$ naturally induces a normal derivative
\[
(Nf)_p:N_p\mathcal F\longrightarrow N_{f(p)}\mathcal F'.
\]
Let $\zeta$ and $\zeta'$ be transverse vector fields on $M$ and $M'$,
respectively. We say that $\zeta$ and $\zeta'$ are $f$-\textbf{related} if
\[
\zeta'_{f(p)}=(Nf)_p\zeta_p,\qquad \forall\,p\in M.
\]

\begin{definition}
Suppose that there are two transverse Lie algebra actions
$a:\mathfrak g\to\mathfrak X(M,\mathcal F)$ and
$a':\mathfrak g\to\mathfrak X(M',\mathcal F')$. We say that a foliate map
$f:M\to M'$ is $\mathfrak g$-\textbf{equivariant} if, for every
$\xi\in\mathfrak g$, the transverse vector fields $a(\xi)$ and $a'(\xi)$
are $f$-related.
\end{definition}

\begin{definition}
Let $(M,\mathcal F)$ be a foliated manifold equipped with a transverse
$\mathfrak g$-action. We say that the local flow $\{\varphi_t\}$ of a
foliate vector field is $\mathfrak g$-\textbf{equivariant} if $\varphi_t$ is
$\mathfrak g$-equivariant for all $t$.
\end{definition}

Assume that there is a transverse Lie algebra action
$a:\mathfrak g\to\mathfrak X(M,\mathcal F)$. For a Lie subalgebra
$\mathfrak h$ of $\mathfrak g$, define
\[
M^{\mathfrak h}
=
\{p\in M\,\vert\, a(\xi)(p)=0,\ \forall\,\xi\in\mathfrak h\}
\]
to be \textbf{the set of fixed leaves} under the transverse action of
$\mathfrak h$. Moreover, denote by $(\mathfrak h)$ the set of Lie
subalgebras of $\mathfrak g$ that are conjugate to $\mathfrak h$. We say that
$x$ is of \textbf{orbit type} $(\mathfrak h)$ if
$\text{stab}(x,\mathfrak g\ltimes\mathcal F)$ is an element of
$(\mathfrak h)$. Define
\[
M_{(\mathfrak h)}
=
\{x\in M\,\vert\, x\text{ is of orbit type }(\mathfrak h)\}.
\]
\subsection{Transverse isometric Lie algebra action}

\begin{definition}\label{isometric} Let $\mathfrak{g}$ be a finite-dimensionalLie algebra, and $g$ a transverse Riemannian metric on a foliated manifold $(M, \mathcal{F})$. A transverse isometric $\mathfrak{g}$-action on $(M, \mathcal{F}, g)$ is a Lie algebra homomorphism 
\begin{equation}\label{isometric-action} a: \mathfrak{g}\rightarrow \mathfrak{X}(M, \mathcal{F}, g).\end{equation}

\end{definition}

Assume that there is a transverse isometric Lie algebra $\mathfrak{g}$-action on a Riemannian foliation $(M, \mathcal{F})$. The following results are proved in \cite{LS18}.

\begin{proposition}\label{linearization}(\cite[Prop. 4.2.]{LS18})  Let $x\in M$ and $\mathfrak{h} = \text{stab}(x, \mathfrak{g}\ltimes F )$. Let $\mathcal{F}_T$ be the linear foliation of the tangent space $T = T_x M$ defined by the linear subspace $F = T_x \mathcal{F}$.
\begin{itemize} \item[a)]  The inner product $g_x$ on $T/F$ defines a transverse Riemannian metric on $(T, \mathcal{F}_T)$. The Lie algebra $\mathfrak{h}$ acts transversely on $(T, \mathcal{F}_T )$ by linear infinitesimal isometries. \item[b)] There is an $\mathfrak{h}$-equivariant foliate open embedding $\psi : T\rightarrow M$ with the properties
$\psi(0) = x$ and $T_0\psi = id_T$. \end{itemize}
\end{proposition}

\begin{theorem}(\cite[Thm. 4.7]{LS18})\label{fixed-leave} Suppose that $\mathfrak{h}$ is a Lie subalgebra of $\mathfrak{g}$ with normalizer $\mathfrak{n}=\mathfrak{n}_{\mathfrak{g}}(\mathfrak{h})$. Then the fixed-leaf set $M^{\mathfrak{h}}$ is a $\mathcal{F}$-saturated closed submanifold of $M$ that is invariant under the transverse action of $\mathfrak{n}$. 
\end{theorem}


In this paper we will focus on the case that $\mathfrak{g}$ is abelian. The following result is an immediate consequence of \cite[Thm. 4.10]{LS18}, and will play an important role in Section \ref{pre-sympl}.

\begin{corollary}\label{finite-isotropy} Suppose that there is a transverse isometric action of a finite-dimensionalabelian Lie algebra $\mathfrak{g}$ on a Riemannian foliation $(M, \mathcal{F})$, and that $M/\overline{\mathcal F}$ is compact. Then the set $\{ \text{stab}(x, \mathfrak{g}\ltimes \mathcal{F})\,\vert\, x\in M\}$ is finite.

\end{corollary}

Let $X$ be a $(\mathfrak{g} \ltimes \mathcal{F})$ -invariant embedded submanifold of $M$, and let $\mathcal{F}_X = \mathcal{F}\vert_X$ be the restriction of the foliation to $X$. Then the normal bundle
\[NX =TM\vert_X /TX\cong N\mathcal{F}/N\mathcal{F}_X\]
is a foliated vector bundle over $(X, \mathcal{F}_X )$ and is equipped with a natural transverse $\mathfrak{g}$-action with the property that the bundle projection $N X \rightarrow X$ is $\mathfrak{g}$-equivariant. A $(\mathfrak{g}\ltimes \mathcal{F})$-invariant tubular neighborhood of X is a
 $\mathfrak{g}$-equivariant foliate embedding $\phi : NX \hookrightarrow M$ with the following properties: the image $\phi(NX)$ is a 
$(\mathfrak{g}\ltimes \mathcal{F})$-invariant open subset of $M$; $\phi\vert_X =id_X$; and $T_x \phi = id_{N_xX}$ for all $x \in X$. (Here we identify $X$ with the zero section of $NX$ and the normal bundle of $X$ in $NX$ with $NX$.)

\begin{proposition}\cite[Prop. 3.3.4]{LS18} Let $(M, \mathcal{F})$ be a Riemannian foliation equipped with a transverse isometric Lie algebra $\mathfrak{g}$-action. Every $(\mathfrak{g}\ltimes \mathcal{F})$-invariant closed embedded submanifold $X$ of M has a $(\mathfrak{g}\ltimes \mathcal{F})$-invariant tubular neighborhood. For every pair of $(\mathfrak{g}\ltimes \mathcal{F})$-invariant tubular neighborhoods $\phi_0, \phi_1: NX \rightarrow M$ there exists a $\mathfrak{g}$-equivariant foliate isotopy from $\phi_0$ to $\phi_1$.
\end{proposition}

\section{Equivariant basic cohomology and Thom isomorphisms}\label{eq-cohomology}

Throughout this section, we assume that $(M,\mathcal F)$ is a foliated
manifold equipped with a transverse Lie algebra action
$a:\mathfrak g\to\mathfrak X(M,\mathcal F)$. For
$\alpha\in\Omega(M,\mathcal F)$, define
\begin{equation}\label{t-action}
\iota(\xi)\alpha=\iota(\xi_M)\alpha,
\qquad
\mathcal L(\xi)\alpha=\mathcal L(\xi_M)\alpha,
\end{equation}
where $\xi_M$ is a foliate vector field representing the transverse vector
field $a(\xi)$. Since $\alpha$ is basic, the contraction and Lie derivative
operators defined above do not depend on the choice of representative
$\xi_M$.

Goertsches and T\"{o}ben \cite[Proposition~3.12]{GT18} observed that these
operators satisfy the usual rules of Cartan's differential calculus, for
example
\[
[\mathcal L(\xi),\mathcal L(\eta)]=\mathcal L([\xi,\eta]).
\]
Equivalently, a transverse $\mathfrak g$-action equips the basic de Rham
complex $\Omega(M,\mathcal F)$ with the structure of a
$\mathfrak g^{\star}$-algebra in the sense of \cite[Chapter~2]{GS99}.
Therefore the $\mathfrak g^{\star}$-algebra $\Omega(M,\mathcal F)$ has a
well-defined Cartan model
\begin{equation}\label{basic-cartan}
\Omega_{\mathfrak g}(M,\mathcal F):=
[S\mathfrak g^*\otimes \Omega(M,\mathcal F)]^{\mathfrak g}.
\end{equation}
An element of $\Omega_{\mathfrak g}(M,\mathcal F)$ may be naturally
identified with an equivariant polynomial map from $\mathfrak g$ to
$\Omega(M,\mathcal F)$, and is called an \textbf{equivariant basic
differential form}.

The equivariant basic Cartan complex has a bigrading given by
\[
\Omega_{\mathfrak g}^{i,j}(M,\mathcal F)=
[S^i\mathfrak g^*\otimes \Omega^{j-i}(M,\mathcal F)]^{\mathfrak g}.
\]
It is equipped with the vertical differential $1\otimes d$, which we
abbreviate to $d$, and the horizontal differential $d'$, defined by
\[
(d'\alpha)(\xi)=-\iota(\xi)\alpha(\xi),
\qquad
\forall\,\xi\in\mathfrak g.
\]
As a single complex, $\Omega_{\mathfrak g}(M,\mathcal F)$ has grading
\[
\Omega_{\mathfrak g}^k(M,\mathcal F)=
\bigoplus_{i+j=k}\Omega_{\mathfrak g}^{i,j}(M,\mathcal F),
\]
and total differential $d_{\mathfrak g}=d+d'$, called the equivariant
exterior differential. The \textbf{equivariant basic de Rham cohomology}
$H_{\mathfrak g}(M,\mathcal F)$ of the transverse $\mathfrak g$-action on
$(M,\mathcal F)$ is defined to be the total cohomology of the Cartan complex
$(\Omega_{\mathfrak g}(M,\mathcal F),d_{\mathfrak g})$.

Similarly, the transverse $\mathfrak g$-action equips
$\Omega_c(M,\mathcal F)$, the de Rham complex of compactly supported basic
forms, with the structure of a $\mathfrak g^{\star}$-algebra. We denote by
$H_{\mathfrak g,c}(M,\mathcal F)$ the cohomology of the corresponding Cartan
model
\[
(\Omega_{\mathfrak g,c}(M,\mathcal F),d_{\mathfrak g}),
\]
and call it the \textbf{compactly supported equivariant basic cohomology}.

\begin{lemma}(\cite[Lemma 4.4.1]{LS21}) \label{f-homotopy} Let $(M,\mathcal{F})$ and $(M', \mathcal{F}')$ be foliated manifolds equipped with transverse actions of a Lie algebra $\mathfrak{g}$. Let $f : [0, 1] \times M \rightarrow M'$ be a $\mathfrak{g}$-equivariant foliate homotopy. Then the pullback morphisms $f_0$ and $f_1:\Omega(M', \mathcal{F}')\rightarrow \Omega(M,\mathcal{F})$ are homotopic as morphisms of $\mathfrak{g}^{\star}$-algebras. In particular they induce the same morphisms in equivariant basic cohomologies: $f_0^* = f_1^*:H_{\mathfrak{g}}(M',\mathcal{F}')\rightarrow H_{\mathfrak{g}}(M,\mathcal{F})$.\end{lemma}

The same argument as used in the proof of \cite[Lemma 4.4.1]{LS21} gives rise to the following result.

\begin{lemma} \label{c-f-homotopy} Let $(M,\mathcal{F})$ and $(M', \mathcal{F}')$ be foliated manifolds equipped with transverse actions of a Lie algebra $\mathfrak{g}$. Let $f : [0, 1] \times M \rightarrow M'$ be a $\mathfrak{g}$-equivariant foliate homotopy.  Suppose that $f$ is a proper map. Then the pullback morphisms $f_0$ and $f_1:\Omega_c(M', \mathcal{F}')\rightarrow \Omega_c(M,\mathcal{F})$ are homotopic as morphisms of $\mathfrak{g}^{\star}$-algebras. In particular they induce the same morphisms in compactly supported equivariant basic cohomologies: $f_0^* = f_1^*:H_{\mathfrak{g}, c}(M',\mathcal{F}')\rightarrow H_{\mathfrak{g},c}(M,\mathcal{F})$.\end{lemma}

For our purpose it will be convenient to work with the following version of the relative basic cohomology.  Let $(M, \mathcal{F})$ be a foliated manifold, and $X \subset M$  a $\mathcal{F}$-saturated submanifold. Define the complex $\Omega(M, X, \mathcal{F})$ to be the kernel of the pullback map $i^*: \Omega(M, \mathcal{F})\rightarrow \Omega(X, \mathcal{F})$. Clearly,  $\Omega(M, X, \mathcal{F})$ is invariant under the usual de Rham exterior differential $d$.  We define the cohomology of the differential complex $(\Omega(M, X, \mathcal{F}), d)$ to be the \textbf{relative basic cohomology} for the pair $(M, X)$, and will have it denoted by $H(M, X, \mathcal{F})$. Similarly, we define the \textbf{relative equivariant basic cohomology} as follows. 


\begin{definition}\label{relative-eq} Let $(M, \mathcal{F})$ be a foliated manifold equipped with a transverse Lie algebra $\mathfrak{g}$-action, and $X \subset M$ a $\mathfrak{g}$-invariant, $\mathcal{F}$-saturated submanifold. Define the complex $\Omega_{\mathfrak{g}}(M, X, \mathcal{F})$ to be the kernel of the pullback map $i^*: \Omega_{\mathfrak{g}}(M, \mathcal{F})\rightarrow \Omega_{\mathfrak{g}}(X, \mathcal{F})$. Then $(\Omega_{\mathfrak{g}}(M, X, \mathcal{F}), d_{\mathfrak{g}})$ is a differential complex. We will define its cohomology $H_{\mathfrak{g}}(M, X, \mathcal{F})$ to be the relative equivariant basic cohomology for the pair $(M, X)$.
\end{definition}

\begin{definition}
Let $(X, \mathcal{F})$ be a foliated manifold, and let $\pi: E\rightarrow X$ be a vector bundle of rank $k$. We say that $E$ is a \textbf{foliated vector bundle}, if there exists an open cover $\{U_{\alpha}\}$ of $X$, and a family of local trivialization map 
\[ \phi_{\alpha}: U_{\alpha}\times \R^k\rightarrow \pi^{-1}(U_{\alpha}),\] 
such that on the overlap $U_{\alpha}\cap U_{\beta}$, the corresponding transition function 
$g_{\alpha\beta}: U_{\alpha} \cap U_{\beta} \rightarrow GL(k)$ is a basic function with respect to the restricted foliation $\mathcal{F}\vert_{U_{\alpha}\cap U_{\beta}}$. We say $E$ is a \textbf{Riemannian foliated vector bundle}, if on the overlap $U_{\alpha}\cap U_{\beta}$, the corresponding transition function $g_{\alpha\beta}$ takes value in $O(k)\subset GL(k)$. 
\end{definition}
If $\pi:E\rightarrow X$ is a foliated vector bundle over $(X, \mathcal{F})$, then $E$ is naturally equipped with a lifted foliation $\mathcal{F}_E$, whose leaves are transverse to the fibers of $\pi$ and are mapped by $\pi$ to those of $\mathcal{F}$.  If $E$ is a Riemannian foliated vector bundle and if $(X, \mathcal{F})$ carries a transversely Riemannian metric $g$, then there is a natural fiberwise Riemannian metric $h$ on $E$, such that $h+\pi^*g$ is a transversely Riemannian metric on $(E, \mathcal{F})$.

\begin{definition} Let $\pi: E\rightarrow X$ be a  foliated vector bundle over a foliated manifold $(X, \mathcal{F})$, and let $X$ be equipped with a transverse  Lie algebra $\mathfrak{g}$-action. We say that $E$ is a $\mathfrak{g}$-equivariant foliated vector bundle, if there is a transverse Lie algebra $\mathfrak{g}$-action on $(E, \mathcal{F}_E)$ such that the bundle map $\pi$ is $\mathfrak{g}$-equivariant.
\end{definition}

Let $X$ be a foliated manifold equipped with a transverse action of a finite-dimensionalLie algebra $\mathfrak{g}$, and $\pi: E\rightarrow X$ an oriented $\mathfrak{g}$-equivariant foliated vector
bundle of rank $r$ over $X$. We say that $A\subset E$ is vertically compact if the restriction $\pi\vert_A: A \rightarrow M$ is a proper map. We say that a basic equivariant differential form in $\Omega_{\mathfrak{g}}(E, \mathcal{F}_E)$ is vertically compactly supported, if its support is a vertically compact subset of $E$.  We will denote by $\Omega^r_{\mathfrak{g}, cv}(E, \mathcal{F}_E)$ the space of vertically compactly supported equivariant basic differential forms on $E$. An \textbf{equivariant basic Thom form} of $E$ is an $r$-form $\tau_{\mathfrak{g}} \in \Omega^r_{\mathfrak{g}, cv}(E, \mathcal{F}_E)$
which satisfies $\pi_*\tau_{\mathfrak{g}} = 1$ and $d_{\mathfrak{g}}\tau_{\mathfrak{g}} = 0$. An equivariant basic Thom form does not always exist for a foliated vector bundle $E$ over $X$. A counterexample was given in \cite[Sec. 4.7]{LS21}. The following result gives a sufficient condition for the existence of an
equivariant basic Thom form on the normal bundle of an invariant submanifold.

\begin{proposition}\label{existence}
{\cite[Prop.~5.2.1]{LS21}}
Let $(M,\mathcal F,g)$ be a Riemannian foliated manifold equipped with an
isometric transverse action of a Lie algebra $\mathfrak g$. Let $X$ be a
co-orientable $(\mathfrak g\ltimes\mathcal F)$-invariant submanifold of
$M$. Then the normal bundle $NX$ possesses an invariant basic metric
connection, and hence admits an equivariant basic Thom form.
\end{proposition}

When the normal bundle admits an equivariant basic Thom form, one obtains the
following equivariant Thom isomorphism in foliated de Rham theory.

\begin{theorem}[Thom isomorphism]\label{thom-isomorphism}
{\cite[Thm.~4.6.1, Prop.~4.8.1]{LS21}}
Let $(X,\mathcal F)$ be a foliated manifold equipped with a transverse
$\mathfrak g$-action, and let $(E,\mathcal F_E,g_E)$ be an oriented
$\mathfrak g$-equivariant Riemannian foliated vector bundle of rank $r$ over
$X$. Suppose that there exists an equivariant basic Thom form
$\tau_{\mathfrak g}$ on $E$. Then fiber integration
\[
\pi_*:\Omega_{\mathfrak g,cv}(E,\mathcal F_E)[r]
\longrightarrow
\Omega_{\mathfrak g}(X,\mathcal F)
\]
is a homotopy equivalence. A homotopy inverse of $\pi_*$ is the Thom map
\[
\zeta_*:\Omega_{\mathfrak g}(X,\mathcal F)
\longrightarrow
\Omega_{\mathfrak g,cv}(E,\mathcal F_E)[r],
\qquad
\zeta_*(\alpha)=\tau_{\mathfrak g}\wedge \pi^*\alpha .
\]
\end{theorem}

The next result will play an important role in the proof of the foliated
Kirwan surjectivity and injectivity theorems.

\begin{proposition}\label{euler}
{\cite[Prop.~6.3.1]{LS18}}
Let $(M,\mathcal F,g)$ be a complete Riemannian foliation equipped with a
transverse isometric action of an abelian Lie algebra $\mathfrak g$, and let
$X$ be a connected component of the fixed-leaf manifold $M^{\mathfrak g}$.
Let $r$ be the codimension of $X$ in $M$, let $E=NX$ be the normal bundle of
$X$ with the fiberwise metric $g_E$ induced by $g$, and let
$\zeta_X:X\hookrightarrow E$ be the inclusion of $X$ as the zero section.
Let $j_x:\{x\}\to X$ be the inclusion of $x\in X$, and let
$a_x:\mathfrak g\to o(E_x,g_{E,x})$ be the induced $\mathfrak g$-action on
the fiber $E_x$. Then the following statements hold.
\begin{itemize}
\item[1)] The bundle $E$ admits a $\mathfrak g$-invariant almost complex
structure $J$. The weights
$\lambda_1,\lambda_2,\ldots,\lambda_l\in\mathfrak g^*$ of the action $a_x$
with respect to $J_x$ are nonzero and independent of $x\in X$. Moreover,
there is a weight space decomposition
\[
E=E_{\lambda_1}\oplus E_{\lambda_2}\oplus\cdots\oplus E_{\lambda_l}
\]
into $\mathfrak g$-equivariant foliated subbundles. In particular, $E$ is
orientable and its rank is $r=2l$.

\item[2)] Let
$\eta_{\mathfrak g}=\zeta_X^*\tau_{\mathfrak g}
\in \Omega_{\mathfrak g}^r(X,\mathcal F|_X)$
be an equivariant basic Euler form with respect to the orientation given by
$J$, and let $\eta_0$ be the component of $\eta_{\mathfrak g}$ in
$S\mathfrak g^*\otimes \Omega^0(X,\mathcal F|_X)$. Then
\[
\eta_0=\lambda_1\lambda_2\cdots\lambda_l\in S^l\mathfrak g^*.
\]
Reversing the co-orientation of $X$ changes the signs of
$\eta_{\mathfrak g}$ and $\eta_0$.

\item[3)] The basic Euler form $\eta_{\mathfrak g}$ becomes invertible in the
algebra $\Omega_{\mathfrak g}(X,\mathcal F|_X)$ after inverting the weights
$\lambda_1,\lambda_2,\ldots,\lambda_l$.
\end{itemize}
\end{proposition}

It follows from Proposition~\ref{euler} that the equivariant basic Euler
class is not a zero divisor, since the weights
$\lambda_1,\ldots,\lambda_l$ are nonzero. For later use, we note that the
condition $X=M^{\mathfrak g}$ can be weakened. Suppose that $\mathfrak h$ is
a one-dimensional subalgebra of the abelian Lie algebra $\mathfrak g$, and
that $X$ is a connected component of the fixed-leaf manifold for the induced
transverse $\mathfrak h$-action. By Theorem~\ref{fixed-leave}, $X$ is a
$(\mathfrak g\ltimes\mathcal F)$-invariant submanifold of $M$. Therefore the
normal bundle $E=NX$ is a foliated vector bundle with an induced transverse
$\mathfrak g$-action. Proposition~\ref{euler} implies that all the weights
$\lambda_1,\ldots,\lambda_l$ are nonzero. These observations give the
following criterion for the equivariant basic Euler class to be nonzero.

\begin{proposition}\label{invertible}
Let $(M,\mathcal F,g)$ be a complete Riemannian foliation equipped with a
transverse isometric action of an abelian Lie algebra $\mathfrak g$. Let $X$
be a connected component of the fixed-leaf manifold of a one-dimensional Lie
subalgebra of $\mathfrak g$, and let $E=NX$ be the normal bundle of $X$.
Then the equivariant basic Euler class $[\eta_X]$ is not a zero divisor in
$H_{\mathfrak g}(X,\mathcal F|_X)$.
\end{proposition}

\section{Equivariant Morse theory on Riemannian foliations}\label{Eq-Morse}

In this section we extend the classical equivariant Morse--Bott theory to
transverse isometric Lie algebra actions on complete Riemannian foliations.
More precisely, we prove a foliated Morse--Bott lemma, show that the
foliated diffeomorphism type of the sublevel set remains unchanged when no
critical level is crossed, and establish a foliated handle presentation
theorem describing what happens when one passes through a critical level.
These results are the Morse-theoretic foundation for the equivariant basic
cohomology arguments developed later in the paper.

We begin with a simple observation about the critical point set of a basic
function. Suppose that $f$ is a basic function on a foliated manifold
$(M,\mathcal F)$, that $x_0$ is a critical point of $f$, and that $L$ is the
leaf through $x_0$. We claim that the leaf closure $\overline L$ is contained
in $\operatorname{Crit}(f)$. To see this, set
\[
A=\{x\in L\,\vert\, (df)_x=0\},
\]
and equip $L$ with the plaque topology. Clearly, $A$ is closed in $L$.
For any $z\in A$, choose a foliation coordinate chart
$(U,\varphi,x^1,\ldots,x^p,y^1,\ldots,y^q)$ around $z$ such that the vertical
plaque through $z$ is given by
\[
y^1=y^1(z),\ldots,y^q=y^k(z).
\]
On the foliation coordinate neighborhood $U$, the function $f\circ\varphi$
depends only on the transverse coordinates $y^1,\ldots,y^k$. Hence
$(df)_z=0$ implies that $df=0$ on the vertical plaque through $z$. Therefore
$A$ is also open in $L$. Since $x_0\in A$, we have $A=L$, and hence
$L\subset\operatorname{Crit}(f)$. Since $\operatorname{Crit}(f)$ is closed,
it follows that $\overline L\subset\operatorname{Crit}(f)$.

\begin{definition}\label{Morse-df}
A basic function $f$ on a foliated manifold $(M,\mathcal F)$ is said to be
Morse--Bott at a critical value $c$ if each connected component of
$\operatorname{Crit}(f)\cap f^{-1}(c)$ is a closed saturated submanifold of
$M$, and if the Hessian of $f$ is nondegenerate in the directions transverse
to these critical submanifolds. We say that $f$ is Morse--Bott on $M$ if it
is Morse--Bott at every critical value.
\end{definition}

Fix a regular value $a$ of $f$. Since $f$ is $\mathfrak g$-invariant, the
sublevel set $M^a$ is invariant under the transverse $\mathfrak g$-action on
$M$. Keep the notation of Section~\ref{review}, let $P^a=\pi^{-1}(M^a)$,
and let $\mathcal F_{P^a}$ be the restriction of $\mathcal F_P$ to $P^a$.
Apply the Molino construction to the restricted foliation on $M^a$. Denote
the corresponding Molino manifold by $W^a$, so that we have the restricted
Molino diagram
\[
M^a \xleftarrow{\ \pi\ } P^a \xrightarrow{\ \rho\ } W^a .
\]
The transverse $\mathfrak g$-action on $M$ induces a transverse isometric
$\mathfrak g$-action on $M^a$. This action lifts to a transverse isometric
$\mathfrak g$-action on $(P^a,\mathcal F_{P^a})$, and projects to a genuine
isometric Lie algebra action of $\mathfrak g$ on $(W^a,g_{W^a})$. Since
$f$ is $\mathfrak g$-invariant, the boundary
$\partial M^a=f^{-1}(a)$ is also invariant under the transverse
$\mathfrak g$-action; hence the induced Lie algebra action on $W^a$ is
tangent to the boundary. Let $G$ be the simply connected Lie group whose Lie
algebra is $\mathfrak g$. Since $M^a$ is compact, $W^a$ is compact. By the
Palais--Lie theorem for manifolds with boundary, the Lie algebra action of
$\mathfrak g$ on $W^a$ integrates to an isometric $G$-action on $W^a$.

The pullback $\pi^*f:P^a\to\mathbb R$ is an
$\mathcal F_{P^a}$-basic function invariant under the actions of both $K$
and $\mathfrak g$. Hence it descends to a function $f_{W^a}$ on $W^a$.
The induced actions of $G$ and $K$ on $W^a$ determine a natural homomorphism
\[
G\times K\longrightarrow \operatorname{Iso}(W^a,g_{W^a}).
\]
The function $f_{W^a}$ is invariant under the image of this homomorphism.
Since $W^a$ is compact, its isometry group is compact. Let $H^a$ be the
closure of the image of $G\times K$ in $\operatorname{Iso}(W^a,g_{W^a})$.
Then $H^a$ is compact, and clearly $f_{W^a}$ is $H^a$-invariant. A routine
check gives the following result.

\begin{lemma}\label{correspondence}
The correspondence
\[
\Omega^0(M^a,\mathcal F|_{M^a})^{\mathfrak g}
\longrightarrow
\bigl(\Omega^0(W^a)\bigr)^{H^a},
\qquad
f\longmapsto f_{W^a}
\]
is one-to-one. Moreover, $f|_{M^a}$ is a $\mathfrak g$-invariant
Morse--Bott function on $M^a$ if and only if $f_{W^a}$ is an
$H^a$-invariant Morse--Bott function on the Molino manifold $W^a$.
\end{lemma}

The following result is an immediate consequence of Lemma~\ref{correspondence}
and \cite[Lemma~4.9]{W69}.

\begin{proposition}\label{existence-Morse}
Suppose that the dimension of $M^a/\overline{\mathcal F}$ is greater than
zero. Then there exists a $\mathfrak g$-invariant basic function that is
Morse--Bott on $(M^a,\mathcal F|_{M^a})$. Indeed, the space of
$\mathfrak g$-invariant basic Morse--Bott functions is dense in the
$C^{\infty}$-topology in the space of $\mathfrak g$-invariant basic
functions on $(M^a,\mathcal F|_{M^a})$.
\end{proposition}

\begin{remark}
When $\dim(M^a/\overline{\mathcal F})=0$, every connected component of
$M^a$ is a single leaf closure. Without loss of generality, assume that
$M^a$ is connected. Then every leaf of $M^a$ is dense in $M^a$, and every
basic function on $M^a$ is constant.
\end{remark}

Throughout the rest of this section, assume that $g$ is a transverse
Riemannian metric induced by a complete bundle-like metric $g_{TM}$ on
$(M,\mathcal F)$. We will also adopt the following notation. For a Riemannian
vector bundle $E$, we denote by $E(r)$ the open disk bundle of radius $r$
that sits inside $E$. We say that a diffeomorphism $\theta$ from $E(r)$ to an
open neighborhood of the zero section $X$ in $E$ is \textbf{fiberwise} if
$\theta(E(r)\cap E_p)\subset E_p$ for all $p\in X$; we say that $\theta$ is
\textbf{origin preserving} if $\theta(p)=p$ for all $p\in X$. As a first
step, we will extend the following equivariant Morse--Bott lemma to a
foliated setting.

\begin{lemma}(\cite[Lemma~4.1]{W69})\label{Morse--Bott}
Let $\pi:E\rightarrow X$ be a Riemannian $H$-vector bundle for some compact
Lie group $H$, and let $f$ be an $H$-invariant Morse--Bott function on $E$
having $X$, the zero section, as a nondegenerate critical submanifold. For
each $p\in X$, let $E_p^+$ and $E_p^-$ be the positive and negative
subspaces of the quadratic form
$T_p^2f:E_p\times E_p\rightarrow\mathbb R$. Suppose that $X$ is connected
and compact, and that $f(X)=0$. Then
$E^+:=\cup_{p\in X}E_p^+$ and $E^-:=\cup_{p\in X}E_p^-$ are subbundles of
$E$ such that
\begin{equation}\label{splitting}
E=E^+\oplus E^-,
\end{equation}
where the splitting respects the Riemannian metric on $E$; moreover, there
exist $r>0$ and an $H$-equivariant diffeomorphism $\theta$ from $E(r)$ to an
open neighborhood of the zero section in $E$, which is fiberwise and origin
preserving, such that
\begin{equation}\label{c-canonical-form}
f(\theta(x,y))
=
\|x\|^2-\|y\|^2,
\qquad
\forall\, (x,y)\in E_p^+\oplus E_p^-,
\quad p\in X.
\end{equation}
\end{lemma}

We now state and prove the foliated Morse--Bott lemma, one of the main
results stated in the introduction.

\begin{lemma}(\textbf{Foliated Morse--Bott Lemma})\label{f-Morse--Bott}
Let $X$ be a connected, compact and nondegenerate critical submanifold of a
$\mathfrak g$-invariant basic Morse--Bott function $f$. For any $p\in X$,
define $N_p^+X$ and $N_p^-X$ to be the positive and negative subspaces of the
quadratic form $T_p^2f:N_pX\times N_pX\rightarrow\mathbb R$. Then both
$N^+X:=\cup_{p\in X}N_p^+X$ and
$N^-X:=\cup_{p\in X}N_p^-X$ are subbundles of $NX$ such that
\begin{equation}\label{splitting2}
NX=N^+X\oplus N^-X,
\end{equation}
where the splitting respects the Riemannian metric on $NX$ induced by the
given transverse Riemannian metric.

Moreover, if $f(X)=0$, then there exist a
$(\mathfrak g\ltimes\mathcal F)$-equivariant tubular neighborhood
$\phi:NX\rightarrow M$, a positive constant $r>0$, and a
$(\mathfrak g\ltimes\mathcal F)$-equivariant diffeomorphism $\theta$ from
$NX(r)$ to an open neighborhood $U$ of the zero section in $NX$, which is
fiberwise and origin preserving, such that
\begin{equation}\label{canonical-form}
f((\phi\circ\theta)(x,y))
=
\|x\|^2-\|y\|^2,
\qquad
\forall\, (x,y)\in N_pX,\quad p\in X.
\end{equation}
\end{lemma}

\begin{proof}
Choose a regular value $a>0$ of $f$ such that
$X\subset\operatorname{int}M^a$. We use the notation introduced before
Lemma~\ref{correspondence}: $P^a=\pi^{-1}(M^a)$ is equipped with the
restricted lifted foliation $\mathcal F_{P^a}$, and
$M^a\xleftarrow{\ \pi\ }P^a\xrightarrow{\ \rho\ }W^a$ is the corresponding
restricted Molino diagram. Let $X_P=\pi^{-1}(X)$, let
$X_W=\rho(X_P)$, and let $H^a$ be the closure of the image of
$G\times K$ in $\operatorname{Iso}(W^a,g_{W^a})$. Then $X_P$ is a compact
submanifold of $P^a$ invariant under both $K$ and
$\mathfrak g\ltimes\mathcal F_{P^a}$, and $X_W$ is a compact submanifold of
$W^a$ invariant under the induced isometric action of $G\times K$, hence
also under $H^a$. The function $f_P:=\pi^*f$ is
$\mathcal F_{P^a}$-basic and invariant under both $\mathfrak g$ and $K$; it
descends through $\rho$ to a $(G\times K)$-invariant, hence $H^a$-invariant,
function $f_{W^a}$ on $W^a$. Moreover, $f_P$ has $X_P$ as a nondegenerate
critical submanifold, and $f_{W^a}$ has $X_W$ as a nondegenerate critical
submanifold.

Let $\phi_W:NX_W\rightarrow W^a$ be an $H^a$-invariant tubular neighborhood
of $X_W$. By Lemma~\ref{Morse--Bott}, there exist $r>0$ and an
$H^a$-equivariant, fiberwise and origin-preserving diffeomorphism
$\theta_W$ from $NX_W(r)$ to an open neighborhood $U$ of the zero section in
$NX_W$, such that $f_{W^a}(\phi_W\circ\theta_W)$ has the form given in
\eqref{c-canonical-form}. Pulling back $\phi_W:NX_W\rightarrow W^a$ through
$\rho:P^a\rightarrow W^a$, we obtain an embedding
$\phi_P:\rho^*NX_W\cong NX_P\rightarrow P^a$, which is a
$(\mathfrak g\ltimes\mathcal F_{P^a})$-equivariant and $K$-invariant
tubular neighborhood of $X_P$. Hence its quotient by $K$ gives rise to a
$(\mathfrak g\ltimes\mathcal F)$-equivariant embedding
$\phi:NX\rightarrow M$. It is easy to see that $\theta_W$ pulls back to a
fiberwise and origin-preserving diffeomorphism $\theta_P$ from $NX_P(r)$
onto an open neighborhood of the zero section in $NX_P$, which is both
$(\mathfrak g\ltimes\mathcal F_{P^a})$-equivariant and $K$-equivariant. As a
result, $\theta_P$ descends to a fiberwise and origin-preserving
diffeomorphism $\theta$ from $NX(r)$ to an open neighborhood of the zero
section in $NX$, which is $(\mathfrak g\ltimes\mathcal F)$-equivariant.
A routine check shows that \eqref{canonical-form} holds.
\end{proof}

We will also need the following elementary observation on gradients of
invariant basic functions.

\begin{lemma}\label{gradient}
Suppose that $a:\mathfrak g\rightarrow\mathfrak X(M,\mathcal F,g)$ is a
transverse isometric action, that $(M,\mathcal F)$ is endowed with a
bundle-like Riemannian metric $g_{TM}$ which is compatible with $g$, and
that $f$ is a $\mathfrak g$-invariant basic function. Then the gradient
vector field $X=\operatorname{grad}(f)$ with respect to $g_{TM}$ is foliate.
Furthermore, the flow generated by $X$ is $\mathfrak g$-equivariant.
\end{lemma}

\begin{proof}
Suppose that $Y\in C^{\infty}(T\mathcal F)$, that $Z$ is a foliate vector
field everywhere perpendicular to the leaves, and that $g$ is the transverse
Riemannian metric associated to $g_{TM}$. To show that $X$ is foliate, it
suffices to show that $g([Y,X],Z)=0$. Since $f$ is a basic function,
$g_{TM}(X,Y)=df(Y)=0$. This implies that the gradient vector field $X$ is
everywhere perpendicular to the leaves. Since both $X$ and $Z$ are
perpendicular to the leaves, we have $g_{TM}(X,Z)=g(X,Z)$. Also note that
$g_{TM}(X,Z)=df(Z)$ is a basic function. Therefore
\[
\begin{split}
0
&=\mathcal L_Y\bigl(g_{TM}(X,Z)\bigr)
 =\mathcal L_Y\bigl(g(X,Z)\bigr) \\
&=(\mathcal L_Yg)(X,Z)+g(\mathcal L_YX,Z)+g(X,\mathcal L_YZ) \\
&=g(\mathcal L_YX,Z)+g(X,\mathcal L_YZ),
\qquad \text{since } \mathcal L_Yg=0, \\
&=g(\mathcal L_YX,Z),
\qquad \text{since } \mathcal L_YZ\in C^{\infty}(T\mathcal F)=\ker(g), \\
&=g([Y,X],Z).
\end{split}
\]
This proves that $\operatorname{grad}(f)$ is foliate. To prove that the flow
of $X$ is $\mathfrak g$-equivariant, it suffices to show that for every
$v\in\mathfrak g$, if the transverse Killing vector field $a(v)$ is
represented by a foliate vector field $\zeta$, then
$[X,\zeta]\in\mathfrak X(\mathcal F)$. Since $f$ is
$\mathfrak g$-invariant, $\mathcal L_{\zeta}f=0$. Thus, for an arbitrary
foliate vector field $Z$ that is perpendicular to the leaves, we have
\[
\mathcal L_{\zeta}\bigl(df(Z)\bigr)
=
(\mathcal L_{\zeta}df)(Z)+df(\mathcal L_{\zeta}Z)
=
df([\zeta,Z])
=
g(X,[\zeta,Z]).
\]
On the other hand, since $df(Z)=g_{TM}(X,Z)=g(X,Z)$ and since $\zeta$
preserves $g$, we have
\[
\begin{split}
\mathcal L_{\zeta}\bigl(df(Z)\bigr)
&=\mathcal L_{\zeta}\bigl(g(X,Z)\bigr) \\
&=(\mathcal L_{\zeta}g)(X,Z)
  +g(\mathcal L_{\zeta}X,Z)
  +g(X,\mathcal L_{\zeta}Z) \\
&=g([\zeta,X],Z)+g(X,[\zeta,Z]).
\end{split}
\]
This implies that $g([\zeta,X],Z)=0$. Since $Z$ is arbitrary,
$[X,\zeta]\in\mathfrak X(\mathcal F)=\ker(g)$.
\end{proof}

We next show that the foliated diffeomorphism type of the sublevel set does
not change when no critical value is crossed.

\begin{theorem}\label{regular-values}
Suppose that $a<b$ are two real numbers in $f(M)$, and that no critical
values lie in the closed interval $[a,b]$. Then there is a foliation-preserving
$\mathfrak g$-equivariant diffeomorphism from $M^b$ onto $M^a$.
\end{theorem}

\begin{proof}   Suppose that $\epsilon>0$ is sufficiently small 
such that $(a-\epsilon, b+\epsilon)$ consists of regular values of $f$. Then $f^{-1}(a-\epsilon, b+\epsilon)$ is a $(\mathfrak{g}\ltimes \mathcal{F})$-invariant open subset of $M$. It follows easily from the argument given in the proof of \cite[Prop. 3.3.6]{LS18} that there exists
a $(\mathfrak{g}\ltimes \mathcal{F})$-invariant function $\lambda$ on $M$,  which equals
$\dfrac{1}{g( \text{grad}(f), \text{grad}(f) )}$ on $f^{-1}[a, b]$, and which vanishes outside of a compact neighborhood of $f^{-1}[a, b]$ that lies in $f^{-1}(a-\epsilon, b+\epsilon)$.  Define a vector field $X$ on $M$ as follows
 \begin{equation}\label{vector} X_q=\lambda(q) \cdot \text{grad}(f)_q, \,\,\forall\, q\in M.\end{equation}
 Since $\lambda$ is a $\mathfrak{g}$-invariant basic function,  $X$ must be a $\mathfrak{g}$-equivariant foliate vector field that  generates a one parameter subgroup $\{\varphi_t\}$ of $\mathfrak{g}$-equivariant and foliation preserving diffeomorphisms. Define
\[ \Psi: M^b\rightarrow M^a, \,\,\, p\mapsto \varphi_{a-b}(p).\]
It is straightforward to check that $\Psi$ is a  $\mathfrak{g}$-equivariant and foliation preserving diffeomorphism from $M^b$ onto $M^a$.

\end{proof}

\begin{definition}\label{attaching-handle}
Let $C$ and $N$ be foliated manifolds equipped with Riemannian foliations
and transverse isometric Lie algebra $\mathfrak g$-actions, where $N$ is a
foliated manifold with boundary. Let $E^+\to C$ and $E^-\to C$ be
Riemannian foliated $\mathfrak g$-equivariant vector bundles. Denote by
$D^+$ and $D^-$ the disk bundles of $E^+$ and $E^-$, respectively, by
$\mathring D^-$ the open disk bundle of $E^-$, and by $S$ the sphere bundle
of $D^-$. The bundle
\begin{equation}\label{disk-bundle}
D^+\oplus D^-=
\{(v,w)\,\vert\, \|v\|\leq 1,\ \|w\|\leq 1\}
\end{equation}
is called a \textbf{basic handle bundle} of type $(E^+,E^-)$ and index
$\operatorname{rank}(E^-)$.

Let $\tilde N\supset N$ be a foliated manifold with boundary, and let
$H\subset \tilde N$ be a closed subset. We write
$\tilde N=N\cup_{D^+\oplus S}H$, and say that $\tilde N$ arises from $N$ by
$\mathfrak g$-equivariantly attaching a basic handle of type $(E^+,E^-)$ if
the following conditions hold:
\begin{itemize}
\item[1)] there exists a homeomorphism $F:D^+\oplus D^-\to H$;
\item[2)] $\tilde N=N\cup H$;
\item[3)] $F|_{D^+\oplus S}$ is a foliation-preserving
$\mathfrak g$-equivariant diffeomorphism onto $\partial N\cap H$;
\item[4)] $F|_{D^+\oplus\mathring D^-}$ is a foliation-preserving
$\mathfrak g$-equivariant diffeomorphism onto $\tilde N\setminus N$.
\end{itemize}
\end{definition}

 Following the methods developed in \cite{Pa63} and \cite{W69}, we will prove a foliated version of the Morse handle body theorem. We will need the following technical result from \cite[Sec. 11]{Pa63}.

\begin{lemma}\label{tech-lemma} Let $\lambda: \mathbb{R}\rightarrow \mathbb{R}$ be a smooth function that is monotone non-increasing and satisfies $\lambda(t)=1$ if $t\leq \frac{1}{2}$, $\lambda(t)>0$ if $t<1$, and $\lambda(t)=0$ if $t\geq 1$. For
$0\leq s\leq 1$, let $\sigma(s)$ be the unique solution of $\dfrac{\lambda(\sigma)}{1+\sigma}=\dfrac{2}{3}(1-s)$ in $[0, 1]$. Then $\sigma$ is strictly monotone increasing, continuous, $C^{\infty}$ in $[0, 1)$ and $\sigma(0)=\dfrac{1}{2}$, $\sigma(1)=1$. Moreover,
if $\epsilon>0$ and $u^2-v^2\geq -\epsilon$ and $u^2-v^2-\dfrac{3\epsilon}{2}\lambda(\dfrac{u^2}{\epsilon})\leq -\epsilon$, then $u^2\leq \epsilon \sigma(\dfrac{v^2}{\epsilon+u^2})$.

\end{lemma}

Next we explain the argument used in \cite[Sec. 11]{Pa63} can be easily adapted to prove the following result.

\begin{proposition} \label{passing-critical}Let $X$ be a connected, compact and non-degenerate critical submanifold of the $\mathfrak{g}$-invariant basic Morse--Bott function $f$, and let
$\phi: NX\rightarrow M$, $NX(r)$, and $\theta$ be given as in Lemma \ref{f-Morse--Bott}. Consider the splitting $NX=N^+X\oplus N^-X$ as given in (\ref{splitting2}). Define a function $h$ on $B(r):=(\phi\circ\theta)(NX(r))$ as follows.
\begin{equation}\label{locally-modify}(h\circ \phi\circ \theta)(x, y)=(f\circ \phi\circ \theta)(x, y)-\dfrac{3\epsilon}{2}\lambda(\dfrac{\vert\vert x\vert\vert^2}{\epsilon}),\end{equation}

 where $(x, y)\in N_p^+X\oplus N_p^-X$, $p\in X$, and $\lambda:\mathbb{R}\rightarrow \mathbb{R}$ is as in Lemma \ref{tech-lemma}. Then for $3\epsilon <r^2$, 
 \begin{equation} \label{comparison} h^{-1}(-\infty, \epsilon]\cap B(r) = f^{-1}(-\infty, \epsilon]\cap B(r).\end{equation}

 Moreover, $Q:=\{p\in B(r)\,\vert\, h(p)\leq -\epsilon\}$ arises from $Z:=\{p\in B(r)\,\vert\, f(p)\leq -\epsilon\}$
 by  $\mathfrak{g}$-equivariantly attaching a basic handle of type $(N^+X, N^-X)$. 
\end{proposition}

\begin{proof}
First note that
$f^{-1}(-\infty,\epsilon]\cap B(r)\subset
h^{-1}(-\infty,\epsilon]\cap B(r)$ is immediate. Conversely, suppose that
$(\phi\circ\theta)(x,y)\in h^{-1}(-\infty,\epsilon]$ and that
$(h\circ\phi\circ\theta)(x,y)\neq (f\circ\phi\circ\theta)(x,y)$. Then
$\lambda(\|x\|^2/\epsilon)>0$, hence $\|x\|^2<\epsilon$. It follows that
$(f\circ\phi\circ\theta)(x,y)=\|x\|^2-\|y\|^2\leq \|x\|^2<\epsilon$.
This proves \eqref{comparison}.

Now let $D^+\subset N^+X$ and $D^-\subset N^-X$ be the disk bundles as in
Definition~\ref{attaching-handle}, and let
$H=\{p\in B(r)\mid f(p)\geq -\epsilon,\ h(p)\leq -\epsilon\}$. Define
$F:D^+\oplus D^-\rightarrow H$ by
\[
(\theta^{-1}\circ\phi^{-1})\circ F(x,y)
=
\left(
(\epsilon\sigma(\|x\|^2)\|y\|^2+\epsilon)^{1/2}x,\,
(\epsilon\sigma(\|x\|^2))^{1/2}y
\right),
\]
for all $(x,y)\in D^+\oplus D^-$. The condition $3\epsilon<r^2$ ensures
that the right hand side lies in $NX(r)$, since $\|x\|\leq 1$, $\|y\|\leq 1$
and $0\leq \sigma\leq 1$ imply that the squared norm of the image is at most
$3\epsilon$.

Define $G:H\rightarrow D^+\oplus D^-$ by
\[
G(z)=\left(
\dfrac{x}{(\epsilon+\|y\|^2)^{1/2}},
\left(\epsilon\sigma\left(\dfrac{\|x\|^2}{\epsilon+\|y\|^2}\right)\right)^{1/2}y
\right),
\qquad
z=(\phi\circ\theta)(x,y)\in H.
\]
A routine check, identical to the one in \cite[Sec.~11]{Pa63}, shows that
$F$ and $G$ are inverse functions and that $F$ maps
$D^+\oplus S$ onto $H\cap \partial Z$, while mapping
$D^+\oplus \mathring D^-$ diffeomorphically onto $Q\setminus Z$.

It remains only to note that $F$ is foliation-preserving and
$\mathfrak g$-equivariant. Let $R:D^+\oplus D^-\rightarrow NX(r)$ denote
the fiberwise map defined by the right hand side of the formula for
$(\theta^{-1}\circ\phi^{-1})\circ F$. Then $F=\phi\circ\theta\circ R$.
Since $N^+X$ and $N^-X$ are $(\mathfrak g\ltimes\mathcal F)$-invariant
foliated subbundles of $NX$, and since the fiberwise metric is basic and
$\mathfrak g$-invariant, the functions $\|x\|^2$ and $\|y\|^2$ appearing in
the definition of $R$ are basic and $\mathfrak g$-invariant. In a foliated
local trivialization of $N^+X\oplus N^-X$, the map $R$ fixes the base point
and acts fiberwise by scalar multiplication with coefficients depending only
on these basic fiberwise norms. Hence $R$ preserves the lifted foliation and
is $\mathfrak g$-equivariant. Since $\phi$ and $\theta$ are
foliation-preserving and $\mathfrak g$-equivariant by
Lemma~\ref{f-Morse--Bott}, the same is true of $F$. Thus $F$ has all the
properties required in Definition~\ref{attaching-handle}.
\end{proof}

\begin{theorem}  \label{crossing-critical} Let $f$ be a proper $\mathfrak{g}$-invariant basic Morse--Bott function. Suppose that $0$ is a critical value of $f$, that $\epsilon>0$ is so small that $0$ is the only critical value of $f$ in $(-\epsilon, \epsilon)$, and that $X_1,  \cdots, X_s$ are all the connected components of the critical submanifold of $f$ that lie in $f^{-1}(0)$. Then there is a $\mathfrak{g}$-equivariant foliated diffeomorphism from 
$M^{\epsilon}$ to  a $(\mathfrak{g}\ltimes \mathcal{F})$-invariant submanifold $\tilde{M}\subset M$ with boundary that arises from $M^{\epsilon}$ by $\mathfrak{g}$-equivariantly attaching $s$ many basic handles of type $(N^+X_i, N^-X_i)$, $1\leq i\leq s$.  \end{theorem}

\begin{proof}  For each $X_i$, let $\phi_i: NX_i\rightarrow M$ be a $(\mathfrak{g}\ltimes \mathcal{F})$-equivariant tubular neighborhood as given in Lemma \ref{f-Morse--Bott}, and let $B_i(r):=\phi(NX_i(r))$ for $r=5\epsilon$. Without loss of generality, choose $\epsilon$ to be so small that
$B_i(r)\cap B_j(r)=\emptyset$, $ \forall\,  i\neq j$.
Define a function $h$ on $V=f^{-1}(-2\epsilon, \infty)$ as follows. If $z\in B_i(r)$ for some $i$, define $g(z)$ as in (\ref{locally-modify}); if $z\notin \cup_{i=1}^s B_i(r)$, then define $g(z)=f(z)$. We first claim that $g$ is smooth. To see this, it suffices to show that $\forall\, 1\leq i\leq s$, the closure of the set $Z_i:=\{z\in B_i(r)\cap f^{-1}(-2\epsilon, \infty)\,\vert\, f(z)\neq h(z)\}$ is contained in $\mathring{B}_i(r)$, the interior of $B_i(r)$. Write $z=(\phi\circ \theta)(x, y)$, where $(x, y)\in N_p^+X_i\oplus N_p^-X_i$ for some $p\in X_i$. If $z\in Z_i$, then by definition we must have $\vert\vert x\vert\vert^2<\epsilon$. However, $f(z)=\vert\vert x\vert \vert^2-\vert\vert y\vert\vert^2>-2\epsilon$. Hence
$\vert\vert y\vert\vert^2<2\epsilon+\vert\vert x\vert\vert^2<3\epsilon$. It follows that $ \overline{Z}_i\subset B_i(\sqrt{4\epsilon})\subset \mathring{B}_i(r)$. This also proves that the function $g$ is proper.

Next we claim that $h$ has no critical points in
$h^{-1}[-\frac{5\epsilon}{4},\frac{5\epsilon}{4}]$. To simplify notation,
write $\zeta=\|x\|^2$, $\eta=\|y\|^2$, and
$\tilde h=h\circ\phi\circ\theta$. Then
$\tilde h(\zeta,\eta)=\zeta-\eta-\frac{3\epsilon}{2}
\lambda(\zeta/\epsilon)$, and hence
\[
\frac{\partial\tilde h}{\partial\zeta}
=
1-\frac{3}{2}\lambda'(\zeta/\epsilon)>0,
\qquad
\frac{\partial\tilde h}{\partial\eta}=-1<0.
\]
Here we use that $\lambda$ is non-increasing by construction, so $\lambda'\leq 0$. Therefore
$d\tilde h$ can vanish only when both radial differentials
$d\|x\|^2$ and $d\|y\|^2$ vanish, namely along the zero section $X$. But on
$X$ one has $h=-\frac{3\epsilon}{2}$, so no point of $X$ lies in
$h^{-1}[-\frac{5\epsilon}{4},\frac{5\epsilon}{4}]$. Thus $h$ has no critical
points in this region.
Now observe that both $h^{-1}[-\epsilon, \epsilon]$ and $h^{-1}(-\frac{5\epsilon}{4}, \frac{5\epsilon}{4})$ are $(\mathfrak{g}\ltimes \mathcal{F})$-invariant sets. Thus there exists a compactly supported $\mathfrak{g}$-invariant basic function $\rho$ which equals $\frac{1}{g(\text{grad} (h), \text{grad}( h))}$ on $h^{-1}[-\epsilon, \epsilon]$,  and which vanishes outside $h^{-1}(-\frac{5\epsilon}{4}, \frac{5\epsilon}{4})$. Define a vector field 
$X$ on $V$ as follows.
\[ X_p= \begin{cases} &\rho\cdot \text{grad}(h), \,\text{if}\, p\in h^{-1}(-\frac{5\epsilon}{4}, \frac{5\epsilon}{4});\\ & 0, \text{if}\, p\notin h^{-1}(-\frac{5\epsilon}{4}, \frac{5\epsilon}{4}).\end{cases}  \]
Note that since $X$ is supported in a compact subset in $V$, it naturally extends to a foliate vector field on $M$. It is clear that $X$ is a complete foliate vector field on $V$; moreover, its flow $\{\varphi_t\}$ is $\mathfrak{g}$-equivariant. Since $h^{-1}(-\infty, \epsilon]\cap V = f^{-1}(-\infty, \epsilon]\cap V$, $\varphi_{-2\epsilon}$ is a $\mathfrak{g}$-equivariant foliated
diffeomorphism from $M^{\epsilon}$ onto 
\[M^{-\epsilon} \bigcup \{ p\in V\,\vert\, p\in f^{-1} [-\epsilon, \infty) \cap h^{-1}(-\infty, -\epsilon]\}= M^{-\epsilon}\bigcup \left(\cup_{i=1}^s \{p\in B_i(r)\,\vert\, h(p)\leq -\epsilon\}\right).\]
Theorem \ref{crossing-critical} now follows easily from Proposition \ref{passing-critical}.

\end{proof}

\section{Applications to the equivariant basic cohomology }\label{application1}

Throughout this section, unless otherwise stated, $(M,\mathcal F)$ is a
complete Riemannian foliation that admits a transverse isometric
$\mathfrak g$-action, and $f:M\to\mathbb R$ is a proper
$\mathfrak g$-invariant basic Morse--Bott function which attains a minimum
on $M$. To simplify notation, for any saturated submanifold $N$ of $M$, we
will often simply denote by $\mathcal F$ the restricted foliation
$\mathcal F|_N$ on $N$.

\begin{lemma}\label{extension} Suppose that  $a$ is a regular value of $f$. Then any basic form $\alpha\in \Omega(M^a, \mathcal{F})$ can be extended to a basic form on the entire manifold $M$.

\end{lemma} 

\begin{proof} First we show that $\alpha$ can be extended to a basic form on a saturated open neighborhood $
M^{<a+\epsilon}:=\{x\in M\,\vert\, f(x)< a+\epsilon\}$ for some $\epsilon>0$. Since by assumption the boundary of $M^a$ is compact,
there exists $\epsilon>0$, and finitely many foliation coordinate neighborhoods $U_1, \cdots, U_k$, such that 
$M^{<a+\epsilon}\subset M^a\cup(\cup_{i=1}^k U_i)$. Moreover, since $a$ is a regular value of the basic function $f$, we may assume that on each $U_i$ there is a foliation coordinate chart 
\[\varphi_i: U_i\rightarrow \mathbb{R}^p\times \mathbb{R}^{q}, z\mapsto (x^1(z), \cdots, x^p(z), y^1(z), \cdots, y^{q-1}(z), a-f(z)),\]
such that the restriction $\varphi_i: U_i\cap M^a\rightarrow \mathbb{R}^p\times\mathbb{R}^{q-1}\times \mathbb{R}_{\geq 0}$ gives rise to a foliation coordinate chart on $M^a$. In other words, the last transverse coordinate $y^q$ is chosen to be the canonical boundary defining function $a-f$. On each $U_i\cap M^a$, $\alpha$ has an expression of the form
\[ \alpha= \displaystyle \sum_{l_1<\cdots <l_s} f_{l_1\cdots l_s}(x_1, \cdots x_p, y_1, \cdots, y_q) dy^{l_1}\cdots dy^{l_s},\]
such that each coefficient $f_{l_1\cdots l_s}$ is a basic function on $(U_i\cap M^a, \mathcal{F}\vert_{U_i\cap M^a})$.
We note that each coefficient $f_{l_1\cdots l_s}$ can be extended to a smooth basic function on $U_i$. Indeed, a routine check shows that if we extend each $f_{l_1\cdots l_s}$ in the single variable $y^q=a-f$, while keeping the remaining variables fixed, following the extension formula precisely given in \cite{S64}, then the resulting smooth function is also basic. Thus $\alpha$ extends to a basic form on the entire $U_i$. Moreover, following the construction used in \cite{S64}, and using the fact that the normal coordinate $y^q=a-f$ is canonically defined, it is easy to see that the extensions of $\alpha$ on different foliation coordinate neighborhoods can be chosen in such a way that they agree with each other on the overlap $U_i\cap U_j\neq \emptyset$. Therefore $\alpha$ gets extended to a basic form $\tilde{\alpha}$ on $M^{<a+\epsilon}$.

Now apply \cite[Prop. 3.3.6]{LS18} to get a $(\mathfrak{g}\ltimes\mathcal{F})$-invariant function $\lambda$ which equals $1$ on $M^a$ and which vanishes outside $M^{<a+\epsilon}$. Define
\[ E(\alpha)(z)=\begin{cases} \lambda(z) \cdot \tilde{\alpha}(z), &\text{if }\, z\in M^{<a+\epsilon};\\0, &\text{if}\, z\notin M^{<a+\epsilon}.\end{cases}\]
Then $E(\alpha)$ is a basic form on $M$ whose restriction to $M^a$ is $\alpha$.

\end{proof}
We first record the long exact sequence associated with a triple of sublevel
sets.

\begin{lemma}\label{long-exact-pairs} Suppose that $a_1<a_2<a_3$ are three regular values of the $(\mathfrak{g}\ltimes \mathcal{F})$-invariant Morse function $f$. Then there is a long exact sequence of cohomologies
\begin{equation}\label{l-e-s2}\cdots \!\!\rightarrow\!\! H^{  *-1}_{\mathfrak{g}}(M^{a_2}\!, M^{a_1}\!, \mathcal{F})\!\!\rightarrow \!\! H^*_{\mathfrak{g}}(M^{a_3}\!, M^{a_2}\!, \mathcal{F})\!\!\xrightarrow{p^*}\!\! H^*_{\mathfrak{g}}(M^{a_3}\!, M^{a_1}\!, \mathcal{F})\!\!\xrightarrow{i^*}\!\! H^*_{\mathfrak{g}}(M^{a_2}\!, M^{a_1}\!, \mathcal{F})\!\!\rightarrow\! \cdots.\end{equation}
Here $p^*$ is induced by the natural chain map \[p:\Omega_{\mathfrak{g}}(M^{a_3}, M^{a_2}, \mathcal{F})\rightarrow   \Omega_{\mathfrak{g}}(M^{a_3}, M^{a_1}, \mathcal{F}), \alpha \mapsto \alpha,\]  and $i^*$ is the pullback map induced by the inclusion
$i:M^{a_2}\hookrightarrow M^{a_3}$.
\end{lemma}
\begin{proof}
We claim that the following sequence is exact:
\[
0\rightarrow \Omega_{\mathfrak g}(M^{a_3},M^{a_2},\mathcal F)
\xrightarrow{p}
\Omega_{\mathfrak g}(M^{a_3},M^{a_1},\mathcal F)
\xrightarrow{i^*}
\Omega_{\mathfrak g}(M^{a_2},M^{a_1},\mathcal F)
\rightarrow 0.
\]
The only point that requires proof is the surjectivity of $i^*$. Since $f$
is proper and attains a minimum, each sublevel set $M^{a_j}$ is compact. In
particular, the extension result in Lemma~\ref{extension} applies to the
compact manifolds with boundary $M^{a_j}$.

Let
\[
\alpha\in \Omega_{\mathfrak g}(M^{a_2},M^{a_1},\mathcal F)
=
\bigl(S\mathfrak g^*\otimes \Omega(M^{a_2},M^{a_1},\mathcal F)\bigr)^{\mathfrak g}.
\]
By Lemma~\ref{extension}, there exists
$\beta\in S\mathfrak g^*\otimes \Omega(M^{a_3},M^{a_1},\mathcal F)$ such
that $i^*\beta=\alpha$.

If the transverse action of $\mathfrak g$ lifts to a true Lie algebra action
$\mathfrak g\rightarrow \mathfrak R(\mathcal F)$ induced by a
foliation-preserving smooth action of a compact and connected Lie group $G$ on
$(M,\mathcal F)$, then
\[
\bigl(S\mathfrak g^*\otimes \Omega(M^{a_3},M^{a_1},\mathcal F)\bigr)^{\mathfrak g}
=
\bigl(S\mathfrak g^*\otimes \Omega(M^{a_3},M^{a_1},\mathcal F)\bigr)^G.
\]
Averaging $\beta$ over the Haar measure on $G$ yields a relative equivariant
basic form
\[
\beta_1\in
\bigl(S\mathfrak g^*\otimes \Omega(M^{a_3},M^{a_1},\mathcal F)\bigr)^{\mathfrak g}
\]
such that $i^*\beta_1=\alpha$.

In the general case of an isometric transverse Lie algebra $\mathfrak g$-action,
properness of $f$ ensures that $M^{a_3}$ is compact, hence transversely
compact. Therefore the averaging operator constructed in \cite{L26} applies
on $M^{a_3}$. Applying this operator to $\beta$ gives the desired relative
equivariant basic form $\beta_1$ with $i^*\beta_1=\alpha$. This proves the
surjectivity of $i^*$, and hence the exactness of the sequence.
\end{proof}

The next lemma shows that no new equivariant basic cohomology is created
between two regular levels.

\begin{lemma}\label{trivial-cohomology} Suppose that $a<b$ are two regular values of the $(\mathfrak{g}\ltimes \mathcal{F})$-invariant Morse--Bott  function $f$ such that $(a, b)$ contains no critical values of $f$. Let $i: M^a\xhookrightarrow {} M^b$ be the inclusion map. Then
the following results hold.
\begin{itemize}\item[1)] $i^*: H_{\mathfrak{g}}(M^b, \mathcal{F})\rightarrow H_{\mathfrak{g}}(M^a, \mathcal{F})$ is an isomorphism;\item[2)] $H_{\mathfrak{g}}(M^b, M^a, \mathcal{F})=0$;
\item[3)]  For any regular value $c<a$, $i^*: H_{\mathfrak{g}}(M^b, M^c, \mathcal{F})\rightarrow H_{\mathfrak{g}}(M^a, M^c, \mathcal{F})$ is an isomorphism.\end{itemize}
\end{lemma}

\begin{proof}    Let $X$ be the vector field as defined in (\ref{vector}), let 
$\{\varphi_t\}$ be the one parameter subgroup of $\mathfrak{g}$-equivariant and foliation preserving diffeomorphisms, and let 
$\phi_t:=\varphi_{-t(b-a)}$, $\forall\, 0\leq t\leq 1$. Then $\{\phi_t\}_{0\leq t\leq 1}$ is a family of $\mathfrak{g}$-equivariant and foliation preserving diffeomorphisms, which satisfies $\phi_1(M^b)=M^a$. Note that $\{\phi_t\}_{0\leq t\leq 1}$ provides a 
$\mathfrak{g}$-equivariant and foliated homotopy between $\phi_0\vert_{M^b}=\text{id}_{M^b}$ and $i\circ \phi_1\vert_{M^b}=\phi_1\vert{M^b}$.
Thus $\phi^*_1\circ i^*=\text{id}_{H_{\mathfrak{g}}(M^b, \mathcal{F})}$.  Since $\phi^*_1: H^*_{\mathfrak{g}}(M^a, \mathcal{F})\rightarrow H^*_{\mathfrak{g}}(M^b, \mathcal{F})$ is an isomorphism, $i^*$ must also be an isomorphism. This proves the first assertion.  The second assertion follows from the first assertion together with the long exact sequence (\ref{l-e-s2}).  It remains to show the third assertion. Note that we may assume the vector field $X$ vanishes on $M^c$. Then $\phi_t(z)=z$, $\forall\, z\in M^c$. Thus  $\phi_0^*, (i\circ\phi_1)^*: \Omega_c(M^b, M^c, \mathcal{F})\rightarrow \Omega_c(M^a, M^c, \mathcal{F})$ are homotopic as morphisms of $\mathfrak{g}^{\star}$-differential graded algebras. Assertion 3 now follows from Lemma \ref{f-homotopy}.

\end{proof}

We next establish an identification between relative equivariant basic
cohomology and compactly supported equivariant basic cohomology on the
difference of two sublevel sets.

\begin{proposition}\label{c-supported-cohomology}
Let $a<b$ be two regular values of the
$(\mathfrak g\ltimes\mathcal F)$-invariant Morse--Bott function $f$. Then
the chain map
\begin{equation}\label{c-map}
\Omega_{\mathfrak g,c}(M^b\setminus M^a,\mathcal F)
\longrightarrow
\Omega_{\mathfrak g}(M^b,M^a,\mathcal F),
\qquad
\alpha\longmapsto \widetilde\alpha,
\end{equation}
induces an isomorphism on cohomology. Here $\widetilde\alpha$ denotes the
extension of $\alpha$ by zero to $M^b$.
\end{proposition}

\begin{proof}
Set $U=M^b\setminus M^a$. The map \eqref{c-map} is well defined: if
$\alpha\in\Omega_{\mathfrak g,c}(U,\mathcal F)$, then
$\operatorname{supp}(\alpha)$ is compact in $U$, so its extension by zero
$\widetilde\alpha$ is smooth on $M^b$ and vanishes on $M^a$.

We prove injectivity. Suppose $\alpha\in\Omega_{\mathfrak g,c}(U,\mathcal F)$
is closed and $\widetilde\alpha=d_{\mathfrak g}\beta$ for some
$\beta\in\Omega_{\mathfrak g}(M^b,M^a,\mathcal F)$. Choose $\epsilon>0$ so
that $[a,a+\epsilon]$ contains no critical values and
$\widetilde\alpha|_{M^{a+\epsilon}}=0$. Then
$\beta|_{M^{a+\epsilon}}$ is closed in
$\Omega_{\mathfrak g}(M^{a+\epsilon},M^a,\mathcal F)$, whose cohomology
vanishes by Lemma~\ref{trivial-cohomology}. Hence
$\beta|_{M^{a+\epsilon}}=d_{\mathfrak g}\gamma$ for some
$\gamma\in\Omega_{\mathfrak g}(M^{a+\epsilon},M^a,\mathcal F)$.

Choose a $(\mathfrak g\ltimes\mathcal F)$-invariant basic function $\lambda$
with $\lambda=1$ on $M^{a+\epsilon/2}$ and
$\operatorname{supp}(\lambda)\subset M^{<a+\epsilon}$. Extending
$\lambda\gamma$ by zero to $M^b$, set
$\beta_c=\beta-d_{\mathfrak g}(\lambda\gamma)$. Then $\beta_c|_U$ is
compactly supported and $d_{\mathfrak g}(\beta_c|_U)=\alpha$. Thus
$[\alpha]=0$.

For surjectivity, let $\eta\in\Omega_{\mathfrak g}(M^b,M^a,\mathcal F)$ be
closed. Choose $\epsilon$ and $\lambda$ as above. Since
$H_{\mathfrak g}(M^{a+\epsilon},M^a,\mathcal F)=0$, there exists
$\zeta\in\Omega_{\mathfrak g}(M^{a+\epsilon},M^a,\mathcal F)$ such that
$\eta|_{M^{a+\epsilon}}=d_{\mathfrak g}\zeta$. Extending $\lambda\zeta$ by
zero to $M^b$, set $\eta_c=\eta-d_{\mathfrak g}(\lambda\zeta)$. Then
$\eta_c$ is cohomologous to $\eta$ and vanishes near $M^a$. Hence
$\alpha:=\eta_c|_U$ lies in $\Omega_{\mathfrak g,c}(U,\mathcal F)$, and its
extension by zero is $\eta_c$. Thus \eqref{c-map} is surjective on
cohomology.
\end{proof}


\begin{theorem} \label{surjection-injection} In addition to the assumptions made in Theorem \ref{crossing-critical}, assume that $\mathfrak{g}$ is abelian, and that there is a one-dimensionalsubalgebra $\mathfrak{h}$ of $\mathfrak{g}$ such that $\text{Crit}(f)$ coincides with the fixed-leaf set $M^{\mathfrak{h}}=\coprod_{i=1}^s X_i$. Then we have the following short exact sequence of $S\mathfrak{g}^*$-modules.
\begin{equation}\label{s-exact-2}
\displaystyle 0\rightarrow \bigoplus_{i=1}^s H^{*-r_i}_{\mathfrak{g}}(X_i, \mathcal{F})\rightarrow H_{\mathfrak{g}}(M^{\epsilon},\mathcal{F})\rightarrow H_{\mathfrak{g}}(M^{-\epsilon}, \mathcal{F})\rightarrow 0.
\end{equation}
Here $r_i$ is the rank of the negative normal bundle $E^-_i$ over $X_i$. 
Moreover, the morphism $i^*: H_{\mathfrak{g}}(M^{\epsilon},\mathcal{F})\rightarrow \bigoplus_{i=1}^s H_{\mathfrak{g}}(X_i, \mathcal{F})$ induced by the inclusion $i: \coprod X_i\rightarrow M^{\epsilon}$ is injective.
\end{theorem}

\begin{proof}
Set
$A=M^{-\epsilon}\cup\bigl(\cup_{i=1}^s
\{p\in B_i(r)\mid g(p)\leq \epsilon\}\bigr)$.
It follows from the proof of Theorem~\ref{crossing-critical} that the
inclusion map $i:A\hookrightarrow M^\epsilon$ induces an isomorphism
$i^*:H^*_{\mathfrak g}(M^\epsilon,\mathcal F)\cong
H^*_{\mathfrak g}(A,\mathcal F)$. Since $M^{-\epsilon}\subset A$ and
$i$ restricts to the identity on $M^{-\epsilon}$, a standard five-lemma
argument applied to the long exact sequences of the pairs
$(M^\epsilon,M^{-\epsilon})$ and $(A,M^{-\epsilon})$ gives
$H^*_{\mathfrak g}(M^\epsilon,M^{-\epsilon},\mathcal F)\cong
H^*_{\mathfrak g}(A,M^{-\epsilon},\mathcal F)$.

By Proposition~\ref{c-supported-cohomology} and
Theorem~\ref{crossing-critical}, we have
\[
H^*_{\mathfrak g}(A,M^{-\epsilon},\mathcal F)\cong
\bigoplus_{i=1}^s
H^*_{\mathfrak g,c}(D_i^+\oplus \mathring D_i^-,\mathcal F).
\]
For each $i$, the projection
$p_i:D_i^+\oplus\mathring D_i^-\to \mathring D_i^-$,
$p_i(x,y)=y$, and the inclusion
$j_i:\mathring D_i^-\hookrightarrow D_i^+\oplus\mathring D_i^-$,
$j_i(y)=(0,y)$, are $\mathfrak g$-equivariant and foliation-preserving. The
fiberwise homotopy $H_{i,t}(x,y)=((1-t)x,y)$ gives a proper
$\mathfrak g$-equivariant foliated deformation retraction from
$D_i^+\oplus\mathring D_i^-$ onto $\mathring D_i^-$. Hence, by
Lemma~\ref{c-f-homotopy},
$H^*_{\mathfrak g,c}(D_i^+\oplus \mathring D_i^-,\mathcal F)\cong
H^*_{\mathfrak g,c}(\mathring D_i^-,\mathcal F)$.

Via the standard fibrewise radial diffeomorphism
$\mathring D_i^-\cong E_i^-$, the latter group is identified with
$H^*_{\mathfrak g,cv}(E_i^-,\mathcal F)$. Thus
Theorem~\ref{thom-isomorphism}, together with Lemma~\ref{long-exact-pairs},
implies that there is a long exact sequence
\begin{equation}\label{l-exact-2}
\cdots\rightarrow
\bigoplus_{i=1}^s H^{*-r_i}_{\mathfrak g}(X_i,\mathcal F)
\rightarrow H^*_{\mathfrak g}(M^\epsilon,\mathcal F)
\rightarrow H^*_{\mathfrak g}(M^{-\epsilon},\mathcal F)
\rightarrow \cdots .
\end{equation}
In the above exact sequence, the map
$H^{*-r_i}_{\mathfrak g}(X_i,\mathcal F)\rightarrow
H^*_{\mathfrak g}(M^\epsilon,\mathcal F)$ factors as
\[
H^{*-r_i}_{\mathfrak g}(X_i,\mathcal F)
\xrightarrow{\zeta_*}
H^*_{\mathfrak g,c}(\mathring D_i^-,\mathcal F)
\longrightarrow
H^*_{\mathfrak g,c}(M^\epsilon\setminus M^{-\epsilon},\mathcal F)
\longrightarrow
H^*_{\mathfrak g}(M^\epsilon,\mathcal F),
\]
where $\zeta_*$ is the Thom isomorphism in
Theorem~\ref{thom-isomorphism}. The middle map is induced by the inclusion of
the $i$-th handle piece into $M^\epsilon\setminus M^{-\epsilon}$, and the
last map is induced by extending compactly supported equivariant basic forms
on $M^\epsilon\setminus M^{-\epsilon}$ by zero to $M^\epsilon$.

Let $i:X_i\hookrightarrow M^\epsilon$ be the inclusion. Then the composition
$H^{*-r_i}_{\mathfrak g}(X_i,\mathcal F)\rightarrow
H^*_{\mathfrak g}(M^\epsilon,\mathcal F)\xrightarrow{i^*}
H^*_{\mathfrak g}(X_i,\mathcal F)$ is given by
$[\alpha]\mapsto [\eta_{X_i,\mathfrak g}]\wedge[\alpha]$, where
$\eta_{X_i,\mathfrak g}$ is the equivariant basic Euler form of
$E_i^-\rightarrow X_i$. By Proposition~\ref{invertible},
$[\eta_{X_i,\mathfrak g}]$ is not a zero divisor. Therefore both
$H^{*-r_i}_{\mathfrak g}(X_i,\mathcal F)\rightarrow
H^*_{\mathfrak g}(M^\epsilon,\mathcal F)$ and
$i^*:H^*_{\mathfrak g}(M^\epsilon,\mathcal F)\rightarrow
\bigoplus_{i=1}^sH^*_{\mathfrak g}(X_i,\mathcal F)$ are injective. This
finishes the proof of Theorem~\ref{surjection-injection}.
\end{proof}

\begin{theorem}\label{Kirwan-surj}(\textbf{Kirwan Surjectivity})
Consider a transverse isometric action of an abelian Lie algebra
$\mathfrak{g}$ on a complete Riemannian foliation $(M,\mathcal{F})$.
Suppose that $f:M\rightarrow \mathbb{R}$ is a proper
$\mathfrak g$-invariant basic Morse--Bott function which attains a minimum on
$M$, that there is a one-dimensional subalgebra $\mathfrak h$ of
$\mathfrak g$ such that $\operatorname{Crit}(f)$ coincides with the
fixed-leaf set $M^{\mathfrak h}$, and that $0$ is a regular value of $f$.
Then the inclusion map $i:Z:=f^{-1}(0)\xhookrightarrow{} M$ induces a
surjective map
\[
\kappa:H_{\mathfrak g}(M,\mathcal F)
\rightarrow
H_{\mathfrak g}(Z,\mathcal F|_Z)
\]
in cohomology.
\end{theorem}
\begin{proof}
Let $F=f^2$. Since $f$ is proper, basic and $\mathfrak g$-invariant, the
function $F$ is also proper, basic and $\mathfrak g$-invariant. Although $F$
may fail to be Morse--Bott at the critical value $0$, it is a basic
Morse--Bott function at every positive critical value. Since $0$ is a regular
value of $f$, we may choose $\epsilon>0$ such that $F$ has no critical value
in $(0,\epsilon]$.

Choose regular values
\[
\epsilon=b_0<b_1<b_2<\cdots,\qquad b_j\to\infty .
\]
For each $j\geq 1$, the function $F$ is Morse--Bott on
$F^{-1}([b_{j-1},b_j])$. Since $F$ is proper, there are only finitely many
critical values of $F$ in $[b_{j-1},b_j]$. Applying
Theorem~\ref{surjection-injection} successively to these critical levels, we see that
the inclusion map $M^{b_{j-1}}\xhookrightarrow{} M^{b_j}$ induces a
surjection
$
H^*_{\mathfrak g}(M^{b_j},\mathcal F)
\rightarrow
H^*_{\mathfrak g}(M^{b_{j-1}},\mathcal F)$.

Let $[\alpha_0]\in H^*_{\mathfrak g}(M^\epsilon,\mathcal F)$, where
$\alpha_0$ is closed. By the above surjectivity, we can choose closed
equivariant basic forms $\alpha_j\in\Omega_{\mathfrak g}(M^{b_j},\mathcal F)$
inductively such that $\alpha_j|_{M^{b_{j-1}}}$ represents
$[\alpha_{j-1}]$. Thus
\[
\alpha_j|_{M^{b_{j-1}}}-\alpha_{j-1}=d_{\mathfrak g}\gamma_{j-1}\]

for some $\gamma_{j-1}\in\Omega_{\mathfrak g}(M^{b_{j-1}},\mathcal F)$.
By Lemma~\ref{extension}, extend $\gamma_{j-1}$ to an equivariant basic form
$\widetilde\gamma_{j-1}$ on $M^{b_j}$, and replace $\alpha_j$ by
$\alpha_j-d_{\mathfrak g}\widetilde\gamma_{j-1}$. Then
$\alpha_j|_{M^{b_{j-1}}}=\alpha_{j-1}$.

Therefore the forms $\alpha_j$ glue to a closed equivariant basic form
$\alpha$ on $M=\cup_jM^{b_j}$ whose restriction to $M^\epsilon$ is
$\alpha_0$. Hence the inclusion map $i:M^\epsilon\xhookrightarrow{}M$
induces a surjection
$
H^*_{\mathfrak g}(M,\mathcal F)
\rightarrow
H^*_{\mathfrak g}(M^\epsilon,\mathcal F)$.

Finally, note that there is a $\mathfrak g$-equivariant and foliation
preserving deformation retraction from $M^\epsilon$ onto $F^{-1}(0)=Z$. It
follows from Lemma~\ref{f-homotopy} that the inclusion map
$i:Z\xhookrightarrow{}M$ induces a surjection
\[
\kappa:H^*_{\mathfrak g}(M,\mathcal F)
\rightarrow
H^*_{\mathfrak g}(Z,\mathcal F|_Z)
\]
in cohomology.
\end{proof}

\begin{theorem}\label{Kirwan-inj}(\textbf{Kirwan Injectivity})
Consider a transverse isometric action of an abelian Lie algebra
$\mathfrak{g}$ on a complete Riemannian foliation $(M,\mathcal{F})$.
Suppose that $f:M\rightarrow \mathbb{R}$ is a proper
$\mathfrak g$-invariant basic Morse--Bott function which attains a minimum on
$M$, and that there is a one-dimensional subalgebra $\mathfrak h$ of
$\mathfrak g$ such that $\operatorname{Crit}(f)$ coincides with the
fixed-leaf set $X:=M^{\mathfrak h}$. Then the inclusion map
$i:X\xhookrightarrow{} M$ induces an injective map
\[
i^*:H^*_{\mathfrak g}(M,\mathcal F)
\rightarrow
H^*_{\mathfrak g}(X,\mathcal F|_X)
\]
in cohomology.
\end{theorem}

\begin{proof}
Choose regular values
\[
b_1<b_2<\cdots,\qquad b_j\to\infty ,
\]
such that $M=\cup_j M^{b_j}$. Since $f$ is proper, each $M^{b_j}$ is compact.
Moreover, since $f$ is Morse--Bott, only finitely many critical values occur
in each $M^{b_j}$. Hence the finite-level induction using
Theorem~\ref{surjection-injection} shows that the inclusion
$i_j:X\cap M^{b_j}\xhookrightarrow{} M^{b_j}$ induces an injective map
\[
i_j^*:H^*_{\mathfrak g}(M^{b_j},\mathcal F)
\rightarrow
H^*_{\mathfrak g}(X\cap M^{b_j},\mathcal F)
\]
for every $j$.

Let $[\alpha]\in H^*_{\mathfrak g}(M,\mathcal F)$ be such that
$i^*[\alpha]=0$. Then the restriction of $[\alpha]$ to
$H^*_{\mathfrak g}(X\cap M^{b_j},\mathcal F)$ vanishes for every $j$.
By the injectivity of $i_j^*$, the restriction of $[\alpha]$ to
$H^*_{\mathfrak g}(M^{b_j},\mathcal F)$ also vanishes for every $j$.

We now choose primitives compatibly. Since $[\alpha]|_{M^{b_1}}=0$, there
exists $\beta_1\in\Omega_{\mathfrak g}(M^{b_1},\mathcal F)$ such that
$d_{\mathfrak g}\beta_1=\alpha|_{M^{b_1}}$. Suppose that
$\beta_{j-1}$ has been chosen on $M^{b_{j-1}}$. Since
$[\alpha]|_{M^{b_j}}=0$, choose
$\beta'_j\in\Omega_{\mathfrak g}(M^{b_j},\mathcal F)$ such that
$d_{\mathfrak g}\beta'_j=\alpha|_{M^{b_j}}$. Then
$\beta'_j|_{M^{b_{j-1}}}-\beta_{j-1}$ is closed.

By applying Theorem~\ref{surjection-injection} successively to the finitely many
critical levels in $[b_{j-1},b_j]$, the restriction map
\[
H^*_{\mathfrak g}(M^{b_j},\mathcal F)
\rightarrow
H^*_{\mathfrak g}(M^{b_{j-1}},\mathcal F)
\]
is surjective. Hence we may choose a closed form
$\theta_j\in\Omega_{\mathfrak g}(M^{b_j},\mathcal F)$ such that
$\theta_j|_{M^{b_{j-1}}}$ represents the class of
$\beta'_j|_{M^{b_{j-1}}}-\beta_{j-1}$. Thus
\[
\beta'_j|_{M^{b_{j-1}}}-\beta_{j-1}-\theta_j|_{M^{b_{j-1}}}
=d_{\mathfrak g}\gamma_{j-1}
\]
for some $\gamma_{j-1}\in\Omega_{\mathfrak g}(M^{b_{j-1}},\mathcal F)$.
By Lemma~\ref{extension}, extend $\gamma_{j-1}$ to an equivariant basic form
$\widetilde\gamma_{j-1}$ on $M^{b_j}$, and define
\[
\beta_j=\beta'_j-\theta_j-d_{\mathfrak g}\widetilde\gamma_{j-1}.
\]
Then $d_{\mathfrak g}\beta_j=\alpha|_{M^{b_j}}$ and
$\beta_j|_{M^{b_{j-1}}}=\beta_{j-1}$.

Therefore the forms $\beta_j$ glue to an equivariant basic form $\beta$ on
$M=\cup_jM^{b_j}$ satisfying $d_{\mathfrak g}\beta=\alpha$. Hence
$[\alpha]=0$, and so $i^*$ is injective.
\end{proof}

\section{Morse inequalities for Riemannian foliations}\label{Morse-inequalities}

Throughout this section, assume that $(M, \mathcal{F})$ is a Riemannian foliation of codimension $n$ on a compact manifold $M$, and that $f: M\rightarrow \mathbb{R}$ is a basic Morse--Bott function. We note that with the only exceptions of Theorem \ref{surjection-injection} and Theorem \ref{Kirwan-surj}, the foundational results established in Section \ref{Eq-Morse} and Section \ref{application1} apply to the current situation with $\mathfrak{g}=\{0\}$. Indeed, Theorem \ref{surjection-injection} and Theorem \ref{Kirwan-surj} depend on Proposition \ref{invertible} which fail to hold in the non-equivariant case in general. However, in view of Theorem \ref{non-orientable-Thom}, the same argument as given in the proof of Theorem \ref{surjection-injection}, together with Theorem  \ref{non-orientable-Thom}, leads to the following result. 

\begin{proposition}\label{relative-cohomology1} Suppose that $0$ is a critical value of $f$, that $\epsilon >0$  is so small that $0$ is the only critical value of $f$ in $(-\epsilon, \epsilon)$, and that $X_1, \cdots, X_s$ are all the connected components of the critical submanifold of $f$ that lies in $f^{-1}(0)$. Then
\[ H^*(M^{\epsilon}, M^{-\epsilon}, \mathcal{F})\cong\bigoplus_{i=1}^s H^{*-r_{i}}(X_i, L_{X_i}, \mathcal{F}).\]
Here $r_{i}$ is the rank of the negative normal bundle of $X_i$, and $L_{X_i}$ is the orientation line bundle of the negative normal bundle $N^-X_i$.
\end{proposition}

To state the Morse inequalities for Riemannian foliations, we first introduce a few notations to set up the stage. For a connected component of the critical submanifold $X$ of $f$, we will denote by $r_X$ the rank of the negative normal bundle of $X$, and $L_X$ the orientation line bundle of the negative normal bundle $N^-X$. Define
\begin{equation}\label{combinitoral-data} \nu_j=\displaystyle \sum_{X \subset \text{Crit}(f)} \text{dim}\left( H^{j-r_{X}}(X, L_{X}, \mathcal{F})\right),
\end{equation}
were $X$ runs over all connected components of critical submanifolds of $f$.
Let $a<b$ be two regular values of  $f$.  $\forall\, 0\leq j\leq n$, we will denote by $b_{j}(M^b, M^a, \mathcal{F})$ the relative basic Betti number $\text{dim}\left( H^{j}(M^b, M^a, \mathcal{F})\right)$, and by $b_j(M, \mathcal{F})$ the basic Betti number $\text{dim}\left(H^{j}(M, \mathcal{F})\right)$.  As an immediate consequence of Lemma \ref{long-exact-pairs} in the case when $\mathfrak{g}$ is trivial, we see that $b_j$ defines a sub-additive function on the collection of all pairs of $(M^b, M^a)$ described above. This leads to the following Morse inequalities.

\begin{theorem}(\textbf{Morse inequalities})\label{Morse-ineq} 
Suppose that $\mathcal{F}$ is a Riemannian foliation on a compact manifold $M$, and that $f$ is a basic Morse--Bott function. Then the following results hold.
\begin{itemize}

\item[a)] $b_j(M,\mathcal F)\leq \nu_j$ for every $0\leq j\leq n$;

\item[b)]$\displaystyle \sum_{j=0}^n (-1)^j b_j(M, \mathcal{F})=  \sum_{j=0}^n (-1)^j \nu_j$;

\item[c)] For every $0\leq k\leq n$,
\[
\sum_{j=0}^k (-1)^{k-j} b_j(M,\mathcal F)
\leq
\sum_{j=0}^k (-1)^{k-j} \nu_j .
\]

\end{itemize}

\end{theorem}





\section{Hamiltonian actions on presymplectic manifolds }\label{pre-sympl}

In this section we consider Hamiltonian transverse actions of a finite-dimensional
abelian Lie algebra $\mathfrak g$ on presymplectic manifolds. As an application
of the results established in Sections~\ref{Eq-Morse} and~\ref{application1},
when the underlying foliation is also Riemannian and when the transverse Lie
algebra action is isometric, we prove Kirwan surjectivity and injectivity
theorems in this setting.

The Kirwan package proved here is restricted to the abelian case. We expect
that, when a transverse isometric Lie algebra action is induced by a
foliation-preserving action of a compact Lie group, the corresponding
nonabelian Kirwan surjectivity and injectivity results can be derived from the
Morse-theoretic foundations developed in this paper, in close analogy with the
classical symplectic setting. One possible approach is to combine these
foundations with the comparison between a compact Lie group action and the
action of its maximal torus, as developed in \cite{M00}. However, the
nonabelian case involves additional technical issues in the foliated setting,
and we leave its treatment for future work.





\begin{definition}
Let $\mathcal{F}$ be a foliation on a smooth manifold $M$.
We say that $\mathcal{F}$ is a \textbf{transversely symplectic foliation},
if there exists a closed 2-form $\omega$, called a \textbf{transversely symplectic form},
such that for each $x\in M$, the kernel of $\omega_x$ coincides with $T_x\mathcal{F}$.
\end{definition}

\begin{definition} (\cite{LS17})\label{G-ham}
Suppose that $(M, \mathcal{F},\omega)$ is a transversely symplectic foliation, and that there is a transverse action of a finite-dimensionalLie algebra $a: \mathfrak{g}\rightarrow \mathfrak{X}(M, \mathcal{F})$.
We say that the action of $\mathfrak{g}$ is \textbf{Hamiltonian},
if there exists a $\mathfrak{g}$-equivariant smooth map $\mu:M\rightarrow\mathfrak{g}^{*}$,
called the moment map, that satisfies the Hamiltonian equation
\begin{equation}\label{ham-action}
\iota(a(\xi))\omega=d\langle \mu,\xi\rangle,\,\,\,\mathrm{for}\,\,\mathrm{all}\,\,\xi\in\mathfrak{g}.
\end{equation}
Here $\langle\cdot,\cdot\rangle$ denotes the dual pairing between $\mathfrak{g}^*$ and $\mathfrak{g}$.
\end{definition}




\begin{lemma}Suppose that $(M, \mathcal{F})$ is a transversely symplectic foliation, and that there is a Hamiltonian transverse Lie algebra action $\rho: \mathfrak{g}\rightarrow \mathfrak{X}(M, \mathcal{F})$ with a moment map $\mu: M\rightarrow \mathfrak{g}^*$. Then the following statements hold true.
\begin{itemize}
\item[a)] \[\text{Crit}(\mu^{\xi})=\{a(\xi)=0\}.\]
\item[b)] \[ d\mu_p(T_pM)=\left(\text{stab}(p, \mathfrak{g}\ltimes \mathcal{F})\right)^{\perp}, \,\forall\, p\in M.\]
Here given a subspace $\mathfrak{h}\subseteq \mathfrak{g}$, $\mathfrak{h}^{\perp}\subseteq \mathfrak{g}^*$ denotes the annihilator of $\mathfrak{h}$.

 \end{itemize}
\end{lemma}
\begin{proof} a) is an immediate consequence from the Hamiltonian equation (\ref{ham-action}). We need only to show b). Suppose that $\xi\in \text{stab}(p, \mathfrak{g}\ltimes \mathcal{F})$.
Then by definition, $a(\xi)(p)=0$. So a) implies that $p$ must be a critical point of
$\mu^{\xi}$.  Thus $\forall\, X\in T_p$, 
\[\langle  d\mu_p(X), \xi\rangle = \langle d \mu^{\xi}, X \rangle_p=0.\]
This shows that $\text{stab}(p, \mathfrak{g}\ltimes \mathcal{F})\subseteq \left(d\mu (T_pM)\right)^{\perp}$. 
Conversely, if $\xi\in \left(d\mu (T_pM)\right)^{\perp}$, then $\langle d \mu^{\xi}, X\rangle_p=\langle  d\mu_p(X), \xi\rangle=0$, $\forall\, X\in T_pM$.
It follows that $p\in \text{Crit}(\mu^{\xi})$. Thus $a(\xi)(p)=0$, which implies that $\xi \in \text{stab}(p, \mathfrak{g}\ltimes \mathcal{F})$.

\end{proof}

\begin{theorem}\label{Morse-function} Suppose that there is a transverse isometric action of a one-dimensionalLie algebra $\mathfrak{g}$ on a transversely symplectic foliation $(M, \mathcal{F}, \omega)$with a given transverse Riemannian metric $g$, and that the action is Hamiltonian with a moment map $f: M\rightarrow \mathfrak{g}^*\cong \mathbb{R}$.  Then
$f$ is a Morse--Bott function of even index whose critical submanifolds coincides with the fixed-leaf set $M^{\mathfrak{g}}$. 
\end{theorem}

\begin{proof}  We compute $M^{\mathfrak{g}}$ in the chart at a point $x\in M^{\mathfrak{g}}$ given by Proposition \ref{linearization}. The tangent space 
$T = T_xM$ is a direct sum $T = E \oplus F$, where $E$ is an inner-product space with a linear isometric $\mathfrak{g}$-action and the foliation of $T$ is the linear foliation $\mathcal{F}_T$ defined by $F$. Thus $ T^{\mathfrak{g}} =E^{\mathfrak{g}} \oplus F$ and $T=E_1\oplus E^{\mathfrak{g}}\oplus F$, where $E^{\mathfrak{g}}$ is the $\mathfrak{g}$-fixed subspace of $E$, and where
$E_1$ is the orthogonal complement of $E^{\mathfrak{g}}$ in $E$ with respect to the inner product.

Let $\psi: T\rightarrow M$ be the $\mathfrak{g}$-equivariant foliate open embedding given in Proposition \ref{linearization},  let $\omega_T=\psi^*\omega$, and let $\omega_0$ be the linear transversely symplectic structure on $T$ given by $\omega\vert_{T_xM}$. Then $\omega_T$ and $\omega_0$ are two presymplectic forms on $T$, both of which have $\mathcal{F}_T$ as their null foliations. Clearly, the transverse action of $\mathfrak{g}$ on $T$ is Hamiltonian with respect to both
$\omega_0$ and $\omega_T$. We will denote the corresponding Hamiltonian functions by $f_0$ and $f_T$, respectively. 
Let $\xi$ be a basis vector in $ \mathfrak{g}$, and let $\overline{\xi}_T$ be the transverse vector field induced by the action of $\xi$. It follows easily from the Hamiltonian equation (\ref{ham-action}) that $p$ is a critical point of $f_0$ or $f_T$ if and only if the transverse vector $\overline{\xi}_T(p)=0$. In particular, this implies that the critical point set for both $f_0$ and $f_T$ is $W:=E^{\mathfrak{g}}\oplus F$. To finish the proof of Theorem \ref{Morse-function}, it suffices to show that the Hessian of $f_T$ at the origin $\mathbf{0}$ is non-degenerate in the transverse direction of $E^{\mathfrak{g}}$, and its index is even.

Now observe that the transverse Riemannian metric $g_x$ on $T$ projects to an inner product $\bar{g}$ 
on $E$, while the transverse action of $\mathfrak{g}$ on $T$ projects to a linear action of $\mathfrak{g}$ on $E$ that preserves the inner product $\bar{g}$. Thus there is a Lie algebra homomorphism $\rho: \mathfrak{g}\rightarrow so(E)$. Let $G$ be the simply connected Lie group whose Lie algebra is
$\mathfrak g$. The homomorphism $\rho:\mathfrak g\rightarrow so(E)$
integrates to a Lie group homomorphism
$
\Theta:G\rightarrow SO(E)$. 
Let $H$ be the closure of $\Theta(G)$ in $SO(E)$. Since $\mathfrak g$ is
abelian, $H$ is a compact connected torus. Moreover, $\omega_0|_E$ and
$\omega_T|_E$ are two $H$-invariant symplectic forms on an $H$-invariant
neighborhood of the origin in $E$, and they agree with each other at
$\mathbf 0$. Therefore, by the standard equivariant Darboux theorem,
$\omega_0|_E$ and $\omega_T|_E$ differ from each other on an open
neighborhood of $\mathbf 0$ in $E$ by an $H$-equivariant symplectomorphism.

However, the symplectic form $\omega_0\vert_E$ and the action of $\mathfrak{g}$ on $E$ are linear, with $E^{\mathfrak{g}}$ being the fixed point set. It follows from the standard result in equivariant symplectic geometry the Hessian of $f_0$ at $x$ is non-degenerate in the transverse direction of $E^{\mathfrak{g}}$, and has even index.  It follows that the Hessian of $f_T$ has the same properties.

\end{proof}



Throughout the rest of this section, we assume that there is a Hamiltonian
transverse action of an abelian Lie algebra $\mathfrak g$ on a transversely
symplectic foliation $(M,\mathcal F,\omega)$ with moment map
$\mu:M\rightarrow \mathfrak g^*$. We also assume that $(M,\mathcal F)$ is a
Riemannian foliation with a given transverse Riemannian metric $g$, and that
the action of $\mathfrak g$ is isometric with respect to $g$.

\begin{lemma}\label{restrict-map}
Let $x,y\in M$ satisfy
$\operatorname{stab}(x,\mathfrak g\ltimes\mathcal F)
=\operatorname{stab}(y,\mathfrak g\ltimes\mathcal F)=\mathfrak h$,
where $\mathfrak h$ is a Lie subalgebra of $\mathfrak g$. If $x$ and $y$ lie
in the same connected component of $M^{\mathfrak h}$, then
\[
\mu(x)+\mathfrak h^{\perp}
=
\mu(y)+\mathfrak h^{\perp}.
\]
As a result, if $M/\overline{\mathcal F}$ is compact, then the set of affine
spaces
\begin{equation}\label{finite-set}
\left\{
\mu(x)+\left(\operatorname{stab}(x,\mathfrak g\ltimes\mathcal F)\right)^{\perp}
\,\vert\, x\in M
\right\}
\end{equation}
is finite.
\end{lemma}

\begin{proof} The map $\mu_{\mathfrak{h}}: M \rightarrow \mathfrak{h}^*$ defined by composing $\mu$ with the projection
$\pi_{\mathfrak{h}^*}: \mathfrak{g}^*\rightarrow \mathfrak{h}^*$ is a moment map for the restricted transverse $\mathfrak{h}$-action.  By definition $\mu_{\mathfrak{h}^*}$ restricts to a locally constant function on $M^{\mathfrak{h}}$, and so $\mu_{\mathfrak{h}^*}(x)=\mu_{\mathfrak{h}^*}(y)$. It follows that 
\[\mu(x)-\mu(y) \in\text{ker}(\pi_{\mathfrak{h}^*})=\mathfrak{h}^{\perp}.\] 
Since $M/\overline{\mathcal{F}}$ is compact, by Corollary \ref{finite-isotropy} the set $\{\text{stab}(x, \mathfrak{g}\ltimes \mathcal{F})\,\vert\,x\in M\}$ is finite.  Since $M^{\mathfrak{h}}$ has only finitely many components for every such Lie subalgebra $\mathfrak{h}$ of $\mathfrak{g}$, (\ref{finite-set}) must be a finite set.

\end{proof}

\begin{lemma}\label{regular} Suppose that $M/\overline{\mathcal{F}}$ is compact, and that $\mathbf{0}\in\mathfrak{g}^*$ is a regular value. Consider the hyperplane Grassmannian $Gr_1(\mathfrak{g})$ parametrizes the set of codimension one subspaces $\mathfrak{h} \subset \mathfrak{g}$. Then the set 
\[U:=\{\mathfrak{h}\in Gr_1(\mathfrak{g})\,\vert \,\mathbf{0} \,\text{is a regular value for}\, \mu_{\mathfrak{h}} :=proj_{\mathfrak{h}^*} \circ \mu\} \]
is open and dense in $Gr_1(\mathfrak{g})$.
 \end{lemma}
\begin{proof} In view of Lemma \ref{restrict-map}, the same argument used in the proof of \cite[Lemma B.3]{BL10} applies to the current situation with minimal changes.
\end{proof}

The following elementary reduction will be used below; its proof is a
straightforward verification and is omitted.

\begin{lemma}\label{restrict-form}
Let $(M,\mathcal{F},\omega)$ be a transversely symplectic foliation, and let
$a:\mathfrak{g}\rightarrow \mathfrak{X}(M,\mathcal{F})$ be a Hamiltonian
transverse action of an abelian Lie algebra $\mathfrak{g}$ with moment map
$\mu:M\rightarrow \mathfrak{g}^*$. Let $\mathfrak{h}$ and $\mathfrak{h}_1$
be Lie subalgebras of $\mathfrak{g}$ such that
$\mathfrak{g}=\mathfrak{h}\oplus\mathfrak{h}_1$. Suppose that $\mathcal{F}$
is also a Riemannian foliation equipped with a transverse Riemannian metric
$g$, that the $\mathfrak{g}$-action is isometric with respect to $g$, and
that $0\in\mathfrak{g}^*$ and $0\in\mathfrak{h}^*$ are regular values of
$\mu$ and $\mu_{\mathfrak{h}}=\operatorname{pr}_{\mathfrak{h}^*}\circ\mu$,
respectively. Set $W:=\mu_{\mathfrak{h}}^{-1}(0)$. Then the following results
hold.
\begin{itemize}
\item[1)] the kernel distribution of $\omega|_W$ defines a foliation
$\mathcal{F}_0$ on $W$, and $(W,\mathcal{F}_0,\omega|_W)$ is transversely
symplectic;
\item[2)] the induced transverse action of $\mathfrak{h}_1$ on
$(W,\mathcal{F}_0,\omega|_W)$ is Hamiltonian;
\item[3)] the restriction of $g$ to $N\mathcal{F}_0$, where
$N\mathcal{F}_0$ is identified with a subbundle of
$N\mathcal{F}|_W$, is a transverse Riemannian metric on
$(W,\mathcal{F}_0)$; moreover, the transverse action of $\mathfrak{h}_1$ on
$W$ is isometric with respect to $g|_{N\mathcal{F}_0}$.
\end{itemize}
\end{lemma}

The following reduction-by-stages statement is an immediate consequence of
Lemma~\ref{regular} and Lemma~\ref{restrict-form}. It will allow us to reduce
Kirwan surjectivity for Hamiltonian transverse actions of abelian Lie algebras
of arbitrary finite dimension to the one-dimensional case established above.

\begin{proposition} \label{reduction-by-stage} Consider the transverse isometric action of an abelian Lie algebra $\mathfrak{g}$ on a Riemannian foliation $(M, \mathcal{F})$ that is also transversely symplectic. Suppose that the action of $\mathfrak{g}$ is Hamiltonian with a moment map $\mu : M \rightarrow \mathfrak{g}^*$ for which $0 \in\mathfrak{g}^*$ is a regular value.
Then we may choose a basis $\xi_1, \cdots, \xi_n$ of 
$\mathfrak{g}$ such that the following results hold. \begin{itemize}
 
  \item[1)] $0\in \mathfrak{g}_k^*$ is a regular value for the map $\mu_k=\text{proj}_{\mathfrak{g}^*_k} \circ \mu$, where $\mathfrak{g}_k=\text{span}\{\xi_1,\cdots, \xi_k\}$, $\forall\, 1\leq k\leq n$.  
  
  \item[2)] For every $0\leq k\leq n-1$, the restriction of
$\mu^{\xi_{k+1}}$ to the submanifold
$M_k:=\mu_k^{-1}(0)$ is a basic Morse--Bott function with critical
submanifold equal to the fixed-leaf set of the transverse action of
$\mathfrak h_{k+1}:=\operatorname{span}\{\xi_{k+1}\}$ on
$(M_k,\mathcal F_k)$, where $\mu_0=0$, $M_0=M$, and $\mathcal F_k$ is the
transversely symplectic foliation on $M_k$ defined by
$\ker(\omega|_{M_k})$.

  
  \end{itemize}
\end{proposition}

We now apply the preceding proposition to the full moment map. The idea is to
impose the equations defining $\mu^{-1}(0)$ one component at a time, and to
apply the one-dimensional Kirwan surjectivity theorem at each stage.

\begin{theorem}\label{kirwan-surj-general}(\textbf{Kirwan Surjectivity})
Consider a transverse isometric action of an abelian Lie algebra
$\mathfrak{g}$ on a Riemannian foliation $(M,\mathcal{F})$ with compact
underlying manifold, and suppose that $(M,\mathcal F)$ is also transversely
symplectic. Suppose that the action of $\mathfrak{g}$ is Hamiltonian with
moment map $\mu:M\rightarrow\mathfrak{g}^*$, and that
$0\in\mathfrak{g}^*$ is a regular value of $\mu$. Then the inclusion map
$i:Z:=\mu^{-1}(0)\xhookrightarrow{}M$ induces a surjective map
$\kappa:H^*_{\mathfrak{g}}(M,\mathcal{F})
\rightarrow H^*_{\mathfrak{g}}(Z,\mathcal{F}|_Z)$ in cohomology.
\end{theorem}

\begin{proof} Let $ \mathfrak{g}_j=\text{span}\{\xi_1,\cdots, \xi_j\}$,  let $\mathfrak{h}_j=\text{span}\{\xi_{j+1}, \cdots, \xi_{n}\}$, where $n=\text{dim}(\mathfrak{g})$, and let $\mu_j=\text{proj}_{\mathfrak{g}^*_j} \circ \mu$ be given as in Proposition \ref{reduction-by-stage}. Set $M_j=\mu_j^{-1}(0)$. To simplify notation, we still denote the restricted foliation $\mathcal{F}\vert_{M_j}$ by $\mathcal{F}$ for every $j$. We see that the inclusion map $i: \mu^{-1}(0)\xhookrightarrow{} M$ factors as a sequence of inclusion map
\[ \mu^{-1}(0)\xhookrightarrow{} \mu_{n-1}^{-1}(0) \xhookrightarrow{} \cdots \xhookrightarrow{} \mu_1^{-1}(0) \xhookrightarrow{} M.\]
Therefore the homomorphism $\kappa: H^*_{\mathfrak{g}}(M, \mathcal{F})\rightarrow H^*_{\mathfrak{g}}(\mu^{-1}(0), \mathcal{F})$ factors as the composition of a sequence of morphisms
\[H^*_{\mathfrak{g}}(M, \mathcal{F})\xrightarrow{i^*} H^*_{\mathfrak{g}}(\mu_{1}^{-1}(0), \mathcal{F})\xrightarrow{i^*} \cdots
\xrightarrow{i^*} H^*_{\mathfrak{g}}(\mu^{-1}(0), \mathcal{F}).\]
To prove Theorem \ref{kirwan-surj-general}, it suffices to show by mathematical induction the following statement: for every $1\leq j\leq n$, the composition map \begin{equation}\label{surjection} \kappa_j: \underbrace{i^*\circ\cdots\circ i^*}_{j-\text{fold}}:  
H^*_{\mathfrak{g}}(M, \mathcal{F})\rightarrow H^*_{\mathfrak{g}}(M_j, \mathcal{F})\end{equation} is surjective, where $M_j=\mu^{-1}_j(0)$.

Clearly, when $j=1$, the statement is an immediate consequence of Theorem \ref{Kirwan-surj}. Indeed, since $M$ is compact, the function $\mu^{\xi_1}$ is proper and attains a minimum on $M$. Assuming that the statement is true for $j-1<n$, we consider the case for $j$. By the inductive hypothesis, the map $\kappa_{j-1}:  
H^*_{\mathfrak{g}}(M, \mathcal{F})\rightarrow H^*_{\mathfrak{g}}(M_{j-1}, \mathcal{F})$ is surjective. Let $\mathcal{F}_{j-1}$ be the presymplectic foliation defined by the kernel of $\omega\vert_{\mu_{j-1}^{-1}(0)}$, and let
\[ \Omega(M_{j-1}, \mathcal{F})_{\text{bas}, \mathfrak{g}_{j-1}}=\{ \alpha \in  \Omega(M_{j-1}, \mathcal{F})\,\vert\, \iota(\xi)\alpha=\mathcal{L}(\xi)\alpha=0, \,\forall\, \xi\in \mathfrak{g}_{j-1}\}.  \]
It is clear from definition that
\[  \Omega(M_{j-1}, \mathcal{F})_{\text{bas}, \mathfrak{g}_{j-1}}= \Omega(M_{j-1}, \mathcal{F}_{j-1}).\]
Thus it follows from \cite[Prop. 3.9]{GT18} that $H^*_{\mathfrak{g}}(M_{j-1}, \mathcal{F})=H^*_{\mathfrak{h}_{j-1}}(M_{j-1}, \mathcal{F}_{j-1})$. A similar argument shows that $H^*_{\mathfrak{g}}(M_{j}, \mathcal{F})=H^*_{\mathfrak{h}_{j-1}}(M_{j}, \mathcal{F}_{j-1})$.
Now applying Theorem \ref{Kirwan-surj} to  the transverse isometric action of $\mathfrak{h}_{j-1}$ on the Riemannian foliation $\mathcal{F}_{j-1}$ and the basic Morse--Bott function $\mu^{\xi_j}\vert_{M_{j-1}}$,   we see that the inclusion map 
$i: M_{j}\hookrightarrow M_{j-1}$ induces a surjection $i^*: H^*_{\mathfrak{h}_{j-1}}(M_{j-1}, \mathcal{F}_{j-1})\rightarrow H^*_{\mathfrak{h}_{j-1}}(M_j, \mathcal{F}_{j-1})$. Here we use the fact that $M_{j-1}$ is compact, and hence $\mu^{\xi_j}\vert_{M_{j-1}}$ is proper and attains a minimum on $M_{j-1}$. This implies the surjectivity of (\ref{surjection}) immediately. 
\end{proof}

\begin{theorem}(\textbf{Kirwan Injectivity}) \label{Kirwan-inj-2}
Consider a transverse isometric action of an abelian Lie algebra
$\mathfrak{g}$ on a Riemannian foliation $(M,\mathcal{F})$ that is also
transversely symplectic. Assume that the action of $\mathfrak{g}$ is
Hamiltonian with moment map $\mu:M\rightarrow\mathfrak{g}^*$. Let
$X:=M^{\mathfrak{g}}$ be the fixed-leaf set, let $i:X\hookrightarrow M$ be
the inclusion map, and let
\begin{equation} \label{pullback}
i^*: H_{\mathfrak{g}}(M,\mathcal{F})\rightarrow
H_{\mathfrak{g}}(X,\mathcal{F}\vert_X)
\end{equation}
be the pullback morphism induced by $i$. Suppose that there exists a generic
element $\xi\in\mathfrak{g}$ such that, for
$\mathfrak h:=\text{span}\{\xi\}$, one has $M^{\mathfrak h}=M^{\mathfrak g}$,
and such that the component $\mu^\xi$ is proper and attains a minimum on
$M$. Then the morphism \eqref{pullback} is injective.
\end{theorem}

\begin{proof}
By the genericity of $\xi$, the fixed-leaf set of
$\mathfrak h:=\text{span}\{\xi\}$ agrees with $M^{\mathfrak g}=X$. Moreover,
by Theorem~\ref{Morse-function}, the component $\mu^\xi$ is a basic
Morse--Bott function whose critical set coincides with $X$. Since $\mu^\xi$
is proper and attains a minimum on $M$, Theorem~\ref{Kirwan-inj} applies to
$f:=\mu^\xi$. Hence the inclusion map $i:X\hookrightarrow M$ induces an
injective map
$i^*:H_{\mathfrak g}(M,\mathcal F)\rightarrow
H_{\mathfrak g}(X,\mathcal F\vert_X)$.
\end{proof}

\begin{remark}
If $M$ has compact underlying manifold, then the properness assumption on
$\mu^\xi$ is automatic. Moreover, $\mu^\xi$ attains both a minimum and a
maximum on $M$.
\end{remark}

\section{ Hamiltonian torus actions on symplectic orbifolds}\label{orbifold}

We first briefly review some basic notions on orbifolds to set up the stage.  We follow the more classical and elementary approach via local charts and atlas and borrow from the expositions in \cite{MM03}, which also serves as a very good introduction to orbifolds from the viewpoint of Lie groupoids. 

let $N$ be a topological space. An \textbf{orbifold chart} of dimension $m\geq 0$ on $N$ is a triple $(U, \Gamma, \phi)$, where $U$ is a connected open set in $\mathbb{R}^m$, $\Gamma$ is a finite group of $\text{Diff}(U)$ and $\phi: U\rightarrow N$ is a continuous and open map 
that induces a homeomorphism $U/\Gamma\rightarrow \phi(U)$. Let $(V, H, \psi)$ be another orbifold chart on $N$. An embedding 
$\sigma: (V, H, \psi)\rightarrow (U, \Gamma, \phi)$ between orbifold charts is an embedding $\sigma: V\rightarrow U$ such that 
$\phi\circ \sigma=\psi$. Two orbifold charts $(U, \Gamma, \phi)$ and $(V, H, \psi)$ are said to be \textbf{compatible} with each other, if $\forall\, x\in \phi(U)\cap \psi(V)$, there exist an orbifold chart $(W, G, \rho)$ with $x\in \rho(W)$, and two embeddings between orbifold charts $\sigma: (W, G, \rho)\rightarrow (U, \Gamma, \phi)$ and $\sigma': (W, G, \rho)\rightarrow (V, H, \psi)$.  An orbifold atlas of dimension $m\geq 0$ of a topological space $N$ is a collection of pairwise compatible orbifold charts $\mathfrak{U}=\{(U_i, \Gamma_i, \phi_i)\}_{i\in I}$ of dimension $m\geq 0$ on $N$, such that $N=\bigcup_{i\in I}\phi_i(U_i)$. An \textbf{orbifold} of dimension $m$ is a pair $(N, \mathfrak{U})$, where $N$ is a second countable Hausdorff topological space and $\mathfrak{U}$ is a maximal orbifold atlas of $N$ of dimension $m$. 
A differential form $\alpha$ on an orbifold $(N, \mathfrak{U})$ is a collection of differential forms $\{\alpha_i\in \Omega(U_i)\,\vert\, (U_i, \Gamma_i, \phi_i)\in \mathfrak{U}\}$, which satisfy the following compatibility condition: if $\sigma_{ji}: (U_j,\Gamma_j, \phi_j)\rightarrow
(U_i, \Gamma_i, \phi_i)$ is an embedding between two orbifold charts in $\mathfrak{U}$, then $\sigma_{ji}^*\alpha_i=\alpha_j$. 
We will denote by $\Omega(N)$ the space of differential forms on the orbifold $N$.
 Clearly,  there is a natural exterior operator $d$ which turns $(\Omega(N), d)$ into a differential complex. The orbifold version of the De Rham theorem asserts the cohomology of the differential complex $(\Omega(N), d)$ is isomorphic to the singular cohomology of $N$ with real coefficients, cf. \cite[Theorem 1]{S56} .
 
 Similarly a vector field on an orbifold $(N, \mathfrak{U})$ is a collection of vector fields $\{Z_i\in \Gamma(T U_i)\, \vert\, (U_i, \Gamma_i, \phi_i)\in \mathfrak{U}\}$,  such that if $\sigma_{ji}: (U_j,\Gamma_j, \phi_j)\rightarrow
(U_i, \Gamma_i, \phi_i)$ is an embedding between two orbifold charts in $\mathfrak{U}$, then vector fields $Z_i$ and $Z_j$ are $\sigma_{ji}$-related.  Let $Z$ be a vector field on an orbifold $N$, and $\Omega(N)$ the space of differential forms on $N$.
It is easy to see that the contraction operator $\iota_Z:\Omega^*(N)\rightarrow \Omega^{*-1}(N)$, the Lie derivative operator
$\mathcal{L}_Z:\Omega^*(N)\rightarrow \Omega^*(N)$, and the exterior operator $d: \Omega^*(N)\rightarrow \Omega^{*+1}(N)$ are well defined and satisfy the usual Cartan identities. 

Let $\mathfrak{U}$ and $\mathfrak{V}$ be two orbifold atlases on two orbifolds $N$ and $Q$ respectively.  A continuous map 
$f: N \rightarrow Q$ is said to be a \textbf{smooth orbifold map} if for any $x\in N$ there are orbifold charts $(U, \Gamma, \phi)\in \mathfrak{U}$ around $x$ and $(V, G, \psi)$ around $f(x)$, such that $f(\phi(U))\subset \psi(V)$ and such that  $f\vert_{\phi(U)}:\phi(U)\rightarrow \psi(V)$ gets lifted to a smooth map $\tilde{f}: U\rightarrow V$ with $\psi\circ \tilde{f}=f\circ \phi$. Let $G$ be a Lie group. A smooth action $a$ of $G$ on an orbifold $N$ is a continuous group action of $G$ on $N$, such that $a: G\times N\rightarrow N$ is a smooth orbifold map.

Suppose that there is a smooth action of a Lie group $G$ on an orbifold $N$. Then every element $\xi\in \mathfrak{g}:=\text{Lie}(G)$ generates a vector field $\xi_N$ on $N$.  $\forall\,\alpha\in\Omega(N)$, we define $\iota_{\xi} \alpha:=\iota_{\xi_N}\alpha$ and $
\mathcal{L}_{\xi}\alpha:=\mathcal{L}_{\xi_N}\alpha$. These two operations, together with the exterior operator $d$, equip
  $\Omega(N)$ with the structure of a $\mathfrak{g}^{\star}$-algebra in the sense of \cite[Ch. 2]{GS99}. We define the $G$-equivariant De Rham cohomology of $N$, denoted by $H_{G}(N, \mathbb{R})$, to be the cohomology of the associated Cartan model $(\Omega_{G}(N):=\left(S\mathfrak{g}^*\otimes \Omega(N)\right)^{\mathfrak{g}}, d_{\mathfrak{g}})$.

Suppose that a torus $T$ with Lie algebra $\mathfrak{t}$ acts on a compact symplectic orbifold $(N, \omega)$ in a Hamiltonian fashion with a moment map $\Phi_N: N\rightarrow \mathfrak{t}^*$.  Choose a $T$-invariant Riemannian metric on $N$ that is compatible with the symplectic form $\omega$. This choice equips the tangent bundle $TN$ with the structure of a Hermitian orbifold vector bundle. Let $\pi: M\rightarrow N$ be the unitary frame bundle of $TN$. This is an orbifold principal bundle over $N$ with structure group $U(k)$, where 
$k=\frac{1}{2}\text{dim}(N)$; moreover, the total space $M$ is a smooth manifold on which the structure group $U(k)$ acts locally free, cf. \cite[Sec. 2.4]{MM03}. Thus $M$ is equipped with the structure of a foliation $\mathcal{F}$ whose leaves are orbits of the $U(k)$-action. By assumption, both the Riemannian metric and the symplectic form on $N$ are $T$-invariant. As a result, the $T$-action on $N$ preserves the unitary frames, and so gets lifted to a $T$-action on $M$. Set $\omega_M=\pi^*\omega$ and $\Phi=\pi^*\Phi_N$. Then $(M, \omega_M, T)$ is a presymplectic Hamiltonian $T$-manifold with a moment map $\Phi$. It is clear that the foliation  $\mathcal{F}$ is both presymplectic and Riemannian, and that the leaf space of $\mathcal{F}$ can be naturally identified with the orbifold $N$. In particular, the results established in Section \ref{pre-sympl} apply to the current situation.
 
Now suppose that $0$ is a regular value of $\Phi_N$. Then it is easy to check that $0$ is also a regular value of $\Phi$. Indeed, 
$\Phi^{-1}(0)$ is a saturated closed submanifold of $M$ whose leaf space is precisely
$\Phi^{-1}_N(0)$. Since $(\Phi^{-1}(0), \mathcal{F}\vert_{\Phi^{-1}(0)})$ is a Riemannian foliation with compact leaves, $\Phi^{-1}_N(0)$ must admit the structure of an orbifold as well. Similarly, let $X$ be the fixed point set of the $T$-action on the orbifold $N$. Then 
$X_M:=\pi^{-1}(X)$ is the set of fixed leaves of the lifted $T$-action on $M$. Choose a generic element $\xi\in \mathfrak{t}$, such that for the one dimension Lie subalgebra $\mathfrak{h}:=\text{span}\{ \xi\}$, $M^{\mathfrak{h}}=M^{\mathfrak{t}}=X_M$. Then it follows from 
Theorem \ref{Morse-function} that $X_M$ is a closed saturated submanifold of $M$. Since all the leaves of $\mathcal{F}\vert_{X_M}$ are compact, the leaf space must be an orbifold as well. In other words,  $X$ must admit the structure of an orbifold. 

Note that the quotient map $\pi: M\rightarrow N$ induces a pullback map $\pi^*: \Omega(N)\rightarrow \Omega(M, \mathcal{F})$, which is an isomorphism of $\mathfrak{t}^*$-algebras. Similarly, we have two other isomorphisms of $\mathfrak{t}^*$-algebras
$\pi^*:  \Omega(\Phi^{-1}_N(0))\rightarrow \Omega(\Phi^{-1}(0), \mathcal{F}\vert_{\Phi^{-1}(0)})$ and 
  $\pi^*:\Omega( X) \rightarrow \Omega(X_M, \mathcal{F}\vert_{X_M})$. So we have three isomorphisms of the equivariant cohomologies $\pi^*: H_{T}(N)\cong H_{\mathfrak{t}}(M, \mathcal{F})$, $\pi^*: H_{T}(\Phi^{-1}(0))\cong H_{\mathfrak{t}}(\Phi^{-1}(0), \mathcal{F}\vert_{\Phi^{-1}(0)})$, 
  and $\pi^*: H_T(X)\cong H_{\mathfrak{t}}(X_M, \mathcal{F}\vert_{X_M})$. Moreover, it is easy to check the following two diagrams commute.

\begin{equation}\label{commu1} \begin{tikzcd}H_{T}(N) \arrow[r, "\kappa"] \arrow[d, " \pi^*"] & H_T(\Phi^{-1}(0)) \arrow[d, "\pi^*"] \\ 
H_{\mathfrak{t}}(M, \mathcal{F}) \arrow[r, "\kappa"]  &H_{\mathfrak{t}}(\Phi^{-1}(0), \mathcal{F}\vert_{\Phi^{-1}(0)}),
\end{tikzcd} \begin{tikzcd} H_{T}(N) \arrow[r, "i^*"] \arrow[d, " \pi^*"] & H_T(X) \arrow[d, "\pi^*"] \\ 
H_{\mathfrak{t}}(M, \mathcal{F}) \arrow[r, "i^*"]  &H_{\mathfrak{t}}(X_M, \mathcal{F}\vert_{X_M}).\end{tikzcd} \end{equation}

In the above two diagrams, the horizontal maps are induced by the inclusions. 
In view of (\ref{commu1}), Theorem \ref{kirwan-surj-general} and Theorem \ref{Kirwan-inj-2} lead to the following results.

\begin{theorem}\label{kir-surj-orbifold}  Consider the Hamiltonian action of a compact and connected torus $T$ on a compact symplectic orbifold $(N, \omega)$ with a moment map $\Phi_N: N\rightarrow \mathfrak{t}^*$. Assume that $0$ is a regular value of $\Phi_N$. Then the pullback map
$\kappa: H_{T}(N)\rightarrow H_{T}(\Phi_N^{-1}(0))$ induced by the inclusion $i:\Phi_N^{-1}(0)\hookrightarrow N$ is surjective.

\end{theorem}

\begin{theorem} \label{kir-inj-orbifold} Consider the Hamiltonian action of a compact and connected torus $T$ on a compact symplectic orbifold $(N, \omega)$. Assume that $X$ is the fixed point set. Then the pullback map
$i^*: H_{T}(N)\rightarrow H_T(X)$ is injective.

\end{theorem}

\section{Examples Beyond Existing Morse--Theoretic Frameworks on Foliations}

The present paper has two main contributions. First, for general transverse
isometric Lie algebra actions on Riemannian foliations, it extends the full
machinery of classical equivariant Morse theory to the foliated setting.
Second, it applies this equivariant Morse theory to transverse Hamiltonian
actions on transversely symplectic foliations that are also Riemannian, and
establishes analogues of the Kirwan injectivity and surjectivity theorems for
transverse isometric actions of abelian Lie algebras.

Goertsches and T\"{o}ben \cite{GT18} developed an equivariant basic Morse
theory for the structural Killing algebra of a Killing foliation. Their
method relies on a special comparison isomorphism. If $(M,\mathcal F)$ is a
Killing foliation with structural Killing algebra $\mathfrak a$, then
\cite[Prop.~4.9]{GT18} gives
\begin{equation}\label{commuting-action}
H_{\mathfrak a}^*(M,\mathcal F)
\cong
H_{\mathfrak a\times\mathfrak{so}(q)}^*(P,\mathcal F_P)
\cong
H_{\mathfrak{so}(q)}^*(W),
\end{equation}
where $P$ is the Molino bundle, $\mathcal F_P$ is the lifted foliation on
$P$, and $W$ is the Molino manifold.

The key point is that \eqref{commuting-action} is specific to the structural
Killing algebra action. It is not available for an arbitrary transverse
isometric Lie algebra action, even when the acting Lie algebra is abelian and
even when the underlying foliation is Killing. By contrast, the
Morse-theoretic framework developed in the present paper applies to arbitrary
transverse isometric Lie algebra actions and does not require the foliation
to be Killing.

We now describe two classes of examples that lie outside the scope of the
structural Killing algebra framework.

\subsection{A Non-Killing Transversely Symplectic Foliation with Morse--Bott Moment Map}

We first construct a compact Riemannian foliation which is transversely
symplectic but not Killing, and which admits a transverse Hamiltonian action
of the non-abelian Lie algebra $\mathfrak{su}(n+1)$.

\begin{proposition}\label{non-killing-example}
There exists a foliation $\mathcal F$ on a compact manifold $M$ such that
$(M,\mathcal F)$ is Riemannian, transversely symplectic, and not Killing.
Moreover, $(M,\mathcal F)$ admits a transverse Hamiltonian action of the
non-abelian Lie algebra $\mathfrak{su}(n+1)$. For a generic
$\xi\in\mathfrak{su}(n+1)$, the basic function
\[
\mu^\xi=\langle\mu,\xi\rangle
\]
is Morse--Bott, with critical submanifolds diffeomorphic to
\[
T^2\times_A S^1.
\]
\end{proposition}
\begin{proof}
Let $A\in SL(2,\mathbb Z)$ be a hyperbolic matrix, with eigenvalues
$\lambda>1$ and $\lambda^{-1}<1$. Let $v_u,v_s\in\mathbb R^2$ be
corresponding eigenvectors, so that
\[
A v_u=\lambda v_u,\qquad A v_s=\lambda^{-1}v_s .
\]
Since $A$ preserves the lattice $\mathbb Z^2$, it induces a diffeomorphism,
still denoted by $A$, of the torus $T^2=\mathbb R^2/\mathbb Z^2$.

Let $F=T^2\times\mathbf{CP}^n$, and define
$\varphi=A\times\operatorname{id}_{\mathbf{CP}^n}:F\to F$. Define $M$ to be
the mapping torus
\[
M=F\times_\varphi S^1=(F\times\mathbb R)/\sim,\qquad
(z,t+1)\sim(\varphi(z),t).
\]
\medskip
\noindent\textbf{1. The induced foliation.}
Define a one-dimensional foliation $\mathcal G$ on
$F=T^2\times\mathbf{CP}^n$ by the rank-one subbundle
\[
E_u\oplus 0\subset T(T^2\times\mathbf{CP}^n),
\qquad
E_u:=\operatorname{span}(v_u)\subset T(T^2).
\]
Since $A_*E_u=E_u$ and $\varphi=A\times\operatorname{id}_{\mathbf{CP}^n}$,
the subbundle $E_u\oplus 0\subset TF$ is $\varphi$-invariant. Hence it
descends to a smooth rank-one subbundle $T\mathcal F\subset TM$. Since every
smooth rank-one distribution is integrable, $T\mathcal F$ defines a
one-dimensional foliation $\mathcal F$ on $M$.

\medskip
\noindent\textbf{2. The transverse symplectic form.}

On $\mathbb R^2$, let $\alpha_u,\alpha_s$ be the linear $1$-forms dual to
$v_u,v_s$. On $T^2\times\mathbb R$, define
\[
\theta_2=\lambda^{-t}\alpha_s,\qquad \theta_3=dt.
\]
Since $A^*\alpha_s=\lambda^{-1}\alpha_s$, the forms $\theta_2$ and
$\theta_3$ are invariant under the deck transformation
$(z,t)\mapsto(\varphi(z),t-1)$. Moreover
\[
d\theta_2=(\log\lambda)\,\theta_2\wedge\theta_3,
\qquad
d\theta_3=0.
\]
Thus $\omega_A=\theta_2\wedge\theta_3$ is closed. It also annihilates the
tangent direction $v_u$ of the lifted foliation.

Let $\omega_{FS}$ be the Fubini--Study form on $\mathbf{CP}^n$, and define
$\Omega=\omega_A+\omega_{FS}$ on $F\times\mathbb R$. Then $d\Omega=0$,
$\iota_{v_u}\Omega=0$, and $\Omega$ is invariant under the deck
transformation. Therefore $\Omega$ descends to a basic closed $2$-form on
$(M,\mathcal F)$. It is nondegenerate on the normal bundle, so
$(M,\mathcal F)$ is transversely symplectic.

\medskip
\noindent\textbf{3. The transverse Hamiltonian Lie algebra action.}
Let $SU(n+1)$ act on $\mathbf{CP}^n$ by the standard projective action and
trivially on $T^2$. Since this action commutes with
$\varphi=A\times\operatorname{id}_{\mathbf{CP}^n}$, it descends to a
foliated action on $M$. Differentiating gives a transverse action of
$\mathfrak{su}(n+1)$ on $(M,\mathcal F)$.

Let
\[
\mu_{\mathbf{CP}^n}:\mathbf{CP}^n\to\mathfrak{su}(n+1)^*
\]
be the standard moment map for the Fubini--Study form, and define
\[
\mu:M\to\mathfrak{su}(n+1)^*,
\qquad
\mu([x,p,t])=\mu_{\mathbf{CP}^n}(p).
\]
Then $\mu$ is well-defined and basic. For every
$\xi\in\mathfrak{su}(n+1)$, if $\xi_M$ denotes the induced transverse vector
field and $\mu^\xi=\langle\mu,\xi\rangle$, then
\[
\iota_{\xi_M}\Omega=d_B\mu^\xi .
\]
Thus the induced transverse $\mathfrak{su}(n+1)$-action is Hamiltonian.

\medskip
\noindent\textbf{4. Morse--Bott components of the moment map.}
Choose a generic element $\xi\in\mathfrak{su}(n+1)$ lying in the Lie algebra
of the diagonal maximal torus, and set $\mu^\xi=\langle\mu,\xi\rangle$.
The function $\langle\mu_{\mathbf{CP}^n},\xi\rangle$ on $\mathbf{CP}^n$ is
Morse, with critical points the fixed points $p_0,\dots,p_n$ of the standard
$SU(n+1)$-action, namely the coordinate lines
$p_j=[e_j]\in\mathbf{CP}^n$. Therefore
\[
\operatorname{Crit}(\mu^\xi)=\bigsqcup_{j=0}^n M_j,
\qquad
M_j=(T^2\times\{p_j\})\times_\varphi S^1.
\]
Each $M_j$ is diffeomorphic to the torus bundle $T^2\times_A S^1$. Since the
Hessian of $\langle\mu_{\mathbf{CP}^n},\xi\rangle$ is nondegenerate at each
$p_j$, the Hessian of $\mu^\xi$ is nondegenerate in the normal directions
transverse to each $M_j$. Hence $\mu^\xi$ is Morse--Bott as a basic function.
For $\mathfrak h=\mathbb R\xi$, the transverse fixed-point set is precisely
the same union, so $\operatorname{Crit}(\mu^\xi)=M^{\mathfrak h}$.

\medskip
\noindent\textbf{5. Failure of the Killing property.}
On the covering $F\times\mathbb R$, the foliation is tangent to the
$v_u$-direction. The transverse vector field $v_s$ commutes with the global
transverse coordinate vector fields and is tangent to the leaf closures.
Thus $v_s$ gives a global trivializing section of the Molino sheaf upstairs.

We claim that the Molino sheaf has no nonzero global section on $M$.
Suppose that $\bar v_s$ is such a section. Its lift to the covering
$F\times\mathbb R$ has the form $\widetilde{\bar v}_s=c\,v_s$. Since this is
a section of the Molino sheaf, the coefficient $c$ is locally constant. As
$F\times\mathbb R$ is connected, $c$ is constant.

Because $\bar v_s$ is defined downstairs, its lift is invariant under the
deck transformation $(z,t)\mapsto(\varphi(z),t-1)$. Therefore
\[
c\,v_s=c\,A_*v_s=c\,\lambda^{-1}v_s .
\]
Since $\lambda^{-1}\neq 1$, it follows that $c=0$. Hence $\bar v_s=0$.

Thus the Molino sheaf has no nonzero global section on $M$. In particular,
it is not globally constant, and $\mathcal F$ is not a Killing foliation.
\end{proof}

This example illustrates two related points. First, it shows that the
Morse-theoretic part of the paper is not confined either to Killing
foliations or to abelian transverse actions: the non-Killing Riemannian
foliation constructed above carries a transverse Hamiltonian
$\mathfrak{su}(n+1)$-action, and for generic
$\xi\in\mathfrak{su}(n+1)$ the component $\mu^\xi$ is a basic
Morse--Bott function whose critical set is the fixed-leaf set of the
one-dimensional subalgebra generated by $\xi$. Second, although the Kirwan
surjectivity and injectivity theorems proved in this paper are stated for
abelian transverse actions, restricting this example to a maximal torus, or
to a generic one-dimensional subalgebra of such a torus, gives nontrivial
abelian Hamiltonian transverse actions satisfying the hypotheses used in the
cohomological applications.

\subsection{Symplectic Toric Quasifolds}
The second class of examples comes from symplectic toric quasifolds, introduced
by Prato as analogues of symplectic toric manifolds. Although a symplectic
toric quasifold can be naturally realized as the leaf space of a Killing
foliation, the natural transverse symmetry associated with its construction is
usually larger than the structural Killing algebra. We explain this point
below, and also indicate how the methods developed in the present paper suggest
a natural approach to computing equivariant basic cohomology rings in this
setting.

We first briefly recall Prato's construction of symplectic toric quasifolds; see
\cite{Pra01} for full details. Fix an $m$-dimensional real vector space $V$.
In Prato's construction, the relevant combinatorial data are encoded by a
triple $(\Delta,Q,\{X_1,\ldots,X_d\})$, where $\Delta\subset V^*$ is an
$m$-dimensional simple convex polytope, $Q\subset V$ is a quasi-lattice, and
$X_1,\ldots,X_d\in Q$ are inward normal vectors to the facets of $\Delta$.
Thus
\begin{equation}\label{polytope}
\Delta=\bigcap_{j=1}^d
\bigl\{\mu\in V^*\,\big|\,\langle\mu,X_j\rangle\geq \lambda_j\bigr\}.
\end{equation}
Here simple means that exactly $m$ facets meet at each vertex, and a
quasi-lattice means the $\mathbb Z$-span of finitely many vectors spanning
$V$ over $\mathbb R$.

Consider the standard Hamiltonian action of $T^d$ on $\mathbb C^d$,
\[
(e^{2\pi i\theta_1},\ldots,e^{2\pi i\theta_d})\cdot
(z_1,\ldots,z_d)
=
(e^{2\pi i\theta_1}z_1,\ldots,e^{2\pi i\theta_d}z_d),
\]
with respect to the standard symplectic form
$\omega_0=\frac{\sqrt{-1}}{2}\sum_{j=1}^d dz_j\wedge d\overline z_j$.
One choice of moment map is
\begin{equation}\label{moment-map}
\Phi:\mathbb C^d\longrightarrow(\mathbb R^d)^*,
\qquad
(z_1,\ldots,z_d)\longmapsto
-\frac12\sum_{j=1}^d |z_j|^2e_j^*+\lambda,
\end{equation}
where $\lambda\in(\mathbb R^d)^*$ is constant and
$e_1^*,\ldots,e_d^*$ is the basis dual to the standard basis
$e_1,\ldots,e_d$ of $\mathbb R^d$.

The normal vectors determine a surjective linear map
\begin{equation}\label{proj}
\pi:\mathbb R^d\longrightarrow V,
\qquad
e_j\longmapsto X_j .
\end{equation}
This map induces an epimorphism of quasi-tori $\Pi:T^d\to D:=V/Q$. Denote its
kernel by $N$. By construction, $N$ is an immersed Lie subgroup of $T^d$.
The restricted $N$-action on $\mathbb C^d$ is Hamiltonian, with moment map
$\Phi_N:\mathbb C^d\to\mathfrak n^*$ obtained by composing $\Phi$ with the
natural projection $(\mathbb R^d)^*\to\mathfrak n^*$, where
$\mathfrak n=\ker\pi$. The value $0$ is regular for $\Phi_N$, so
$M:=\Phi_N^{-1}(0)$ is a smooth submanifold of $\mathbb C^d$. Prato defines
the corresponding \emph{symplectic toric quasifold} to be the quotient
\[
X:=M/N=\Phi_N^{-1}(0)/N .
\]

Without loss of generality\footnote{See \cite[Sec. 6.1]{LY19} for a detailed discussion.},
we may assume that $N$ is connected. Since $\Phi_N$ is $T^d$-invariant,
$M=\Phi_N^{-1}(0)$ is $T^d$-invariant. In particular, the $N$-orbits on $M$
determine a foliation $\mathcal F$. The restriction of $\omega_0$ to $M$
induces a transverse symplectic form on $(M,\mathcal F)$, so $\mathcal F$ is
transversely symplectic. The standard Hermitian metric on $\mathbb C^d$
restricts to a Riemannian metric on $M$. Since the $N$-action is by
isometries, this restricted metric is bundle-like for $\mathcal F$ and hence
induces a transverse metric on $(M,\mathcal F)$.

Write $\mathfrak t$ for the Lie algebra of $T^d$ and set
$\mathfrak g:=\mathfrak t/\mathfrak n$. The infinitesimal $T^d$-action gives
a vector field on $M$ for each element of $\mathfrak t$, while elements of
$\mathfrak n$ give vector fields tangent to $\mathcal F$. Therefore it
induces a Lie algebra homomorphism
\[
\mathfrak g=\mathfrak t/\mathfrak n
\longrightarrow
\mathfrak X(M,\mathcal F),
\]
and hence defines a transverse Lie algebra action of $\mathfrak g$ on
$(M,\mathcal F)$. The restriction of $\Phi$ to $M$ determines a basic map
$\mu:M\to\mathfrak g^*$, which is a moment map for this transverse action.
Clearly, $\dim\mathfrak g=\frac12\operatorname{codim}\mathcal F$.

Denote by $H$ the closure of $N$ in $T^d$, and by $\mathfrak h$ its Lie
algebra. It is explained in \cite[Example 4.3]{GT18} that, in this case,
the foliation $\mathcal F$ is Killing with structural Lie algebra
$\mathfrak h/\mathfrak n$. We note that when $N$ is not dense in $T^d$,
$\mathfrak h/\mathfrak n$ is in general a proper Lie subalgebra of
$\mathfrak g$. This shows that even in the present Killing setting, the
transverse Lie algebra action naturally associated with the quasifold
construction need not coincide with the structural Killing algebra action
considered in \cite{GT18}. Thus the equivariant Morse theory developed in the
present paper has broader applicability than the framework developed in
\cite{GT18}.

By Theorem~\ref{Kirwan-inj-2}, the inclusion
$i:M^{\mathfrak g}\hookrightarrow M$ induces an injection
\begin{equation}\label{image}
i^*: H_{\mathfrak g}(M,\mathcal F)
\longrightarrow
H_{\mathfrak g}(M^{\mathfrak g},\mathcal F).
\end{equation}

In the special case when $\Delta$ is a Delzant polytope, $N$ is closed, $G:=T^d/N$ is a compact connected torus, and $X=M/N$ is a symplectic toric manifold. In this situation, it is easy to see that $H_{\mathfrak g}(M,\mathcal F)\cong H_G(X)$, and the map \eqref{image} becomes the usual localization map
\begin{equation}\label{image2}
i^*: H_G(X)\longrightarrow H_G(X^G),
\end{equation}
where $G=T^d/N$. It is well known that the image of \eqref{image2} admits a classical combinatorial description, obtained using either the Chang--Skjelbred lemma or GKM theory. This classical picture, together with the method developed in the present paper, suggests the possibility of developing a GKM theory in the foliated setting and obtaining a combinatorial description of the equivariant basic cohomology ring $H_{\mathfrak g}(M,\mathcal F)$ for symplectic toric quasifolds. To the best of our knowledge, such a computation has not appeared in the literature. We leave this direction for future research.

\appendix
\section{Thom isomorphism for non-orientable foliated vector bundles}\label{non-orientable}

Let $\pi: E\rightarrow M$ be a non-orientable foliated vector bundle of rank $k$ over a foliated manifold $(M, \mathcal{F})$, let $\mathcal{F}_E$ be the lifted foliation on $E$, and let $\mathcal{O}_E$ be the orientation line bundle of $E$. By definition, $\mathcal{O}_E$ is a flat vector bundle. Recall that the twisted De Rham complexes $(\Omega(M, \mathcal{O}_E), d)$ is defined as follows, cf. \cite[Sec. 7]{BT82}.
As a space,  $\Omega^r(M, \mathcal{O}_E):=\Gamma(\wedge^r (T^*M)\otimes \mathcal{O}_E)$. Choose a family of trivialization
$\{(U_{\alpha}, \varphi_{\alpha})\}$ of the orientation line bundle $\mathcal{O}_E$ such that the corresponding transition functions are locally constant. On each $U_{\alpha}$, a twisted differential form $\gamma$  of degree $r$ has an representation $\gamma=\beta\otimes s$, where $\beta$ is a differential forms on $U_{\alpha}$ and $s$ is a nowhere vanishing local section of $\mathcal{O}_E$ over $U_{\alpha}$. The exterior derivative of $\gamma$ is given by the formula
\begin{equation}\label{exterior} d\gamma :=d\beta\otimes s.\end{equation}
It is straightforward to check that $d\gamma$ does not depend on the choice of a local trivialization, and thus defines a global element in 
$\Omega^{r+1}(M, \mathcal{O}_E)$. Similarly, given a vector field $X$ on $M$, define
\[ \mathcal{L}_X\gamma=(\mathcal{L}_X\beta)\otimes s,\,\,\iota_X\gamma= \iota_X\beta \otimes s.\]
One checks easily the above definition does not depend on the choice of a local trivialization, and
extends the usual Lie derivative $\mathcal{L}_X$ and interior product $\iota_X$ to twisted differential forms. 
We say that $\gamma \in \Omega(M, \mathcal{O}_E)$ is a \textbf{twisted basic differential form}, if $\mathcal{L}_X\gamma=\iota_{X}\gamma=0$, $\forall\, X\in \Gamma(\mathcal{F})$, and will denote by $\Omega(M, \mathcal{O}_E, \mathcal{F})$ the space of twisted basic differential forms. It is clear that $\Omega(M, \mathcal{O}_E, \mathcal{F})$ is invariant under the exterior derivative defined in (\ref{exterior}). We call the cohomology of the differential complex $(\Omega(M, \mathcal{O}_E, \mathcal{F}), d)$ the twisted basic cohomology, and denote it by $H(M, \mathcal{O}_E, \mathcal{F})$. Moreover, the same argument as given in the proof \cite[Prop. 7.5]{BT82} shows up to isomorphisms the twisted basic cohomologies are independent of the choice of local trivializations.  Observe that $\pi^*E$ is a vector bundle over $E$, for which the pullback bundle $\pi^*\mathcal{O}_E$ is the orientation line bundle. Applying the above construction to $E$ and $\pi^*\mathcal{O}_E$, we get the De Rham complex of twisted basic differential forms $(\Omega(E, \pi^*\mathcal{O}_E, \mathcal{F}_E), d)$.

Since $E$ is non-orientable, the integration along the fiber operator $\pi_*: \Omega^*_{cv}(E)\rightarrow \Omega(M)$ is no longer well defined. However,  the usual construction in the orientable case can be easily modified to define two integration along the fiber operators in most obvious ways: $\pi_*: \Omega_{cv}(E)\rightarrow \Omega(M,\mathcal{O}_E)$, and $\pi_*: \Omega_{cv}(E, \mathcal{O}_E)\rightarrow \Omega(M)$.  It is easy to check that both operators satisfy properties similar to those stated in \cite[Lemma B.2]{LS21}, and thus induce two operators $\pi_*: \Omega(E, \mathcal{F}_E)\rightarrow \Omega(M, \mathcal{O}_E, \mathcal{F})$ and $\pi_*: \Omega(E, \pi^*\mathcal{O}_E, \mathcal{F}_E)\rightarrow \Omega(M, \mathcal{F})$.  In this context,   a \textbf{ twisted basic Thom form} is a form $\tau \in  \Omega^k_{cv}(E, \pi^*\mathcal{O}_E, \mathcal{F}_E)$ such that $\pi_*(\tau)=1$.
It is easy to see that when the orthogonal frame bundle $P$ of $E$ admits a connection that is  $\mathcal{F}_P$-basic,  a slight modification of the argument given in the proof of \cite[Prop. 4.8.2]{LS21} shows $E$ will admit a twisted basic Thom form. 

Now suppose that $\gamma$ and $\gamma'$ are two twisted basic forms in $ \Omega(E,\pi^*\mathcal{O}_E, \mathcal{F}_E)$. We can define $\gamma \wedge \gamma'\in \Omega(E, \pi^*(\mathcal{O}_E\otimes\mathcal{O}_E), \mathcal{F}_E)$ as follows.
 Set $\gamma=\beta \otimes s$ and $\gamma'=\beta'\otimes s$, where $\beta$ and $\beta'$ are differential forms on $U$ and $s$ is a local section of $\pi^*\mathcal{O}_E$ over $U$. Then on $U$ 
 \[ \gamma \wedge\gamma' =( \beta\wedge\beta')\otimes (s\otimes s).\]
It is straightforward to check that the above definition does not depend on the choice of a local section $s$, and so $\gamma \wedge \gamma'$ is globally defined. Moreover, since $\pi^*(\mathcal{O}_E\otimes \mathcal{O}_E)$ is a trivial line bundle over $E$, $\Omega(E, \pi^*(\mathcal{O}_E\otimes\mathcal{O}_E), \mathcal{F}_E)=\Omega(E, \mathcal{F}_E)$. 
Thus $\gamma\wedge \gamma'$ is a basic differential form on $E$. In particular, we see that if $\tau$ is a twisted basic  Thom form
on $E$, and $\sigma \in \Omega(M, \mathcal{O}_E, \mathcal F)$, then $\tau \wedge \pi^*\sigma \in \Omega_{cv}(E, \mathcal{F}_E)$. 

After these preparations,  the arguments used in \cite{LS21} to establish the Thom isomorphism for an oriented foliated vector bundle can be easily modified to show the following result. . 

\begin{theorem} \label{non-orientable-Thom} Let $(X, \mathcal{F})$ be a foliated manifold, and let $(E,\mathcal F_E,g_E)$ be a possibly non-orientable Riemannian foliated vector bundle of rank $r$ over $X$. Suppose that the orthogonal frame bundle $P$ of $E$ admits a $\mathcal{F}_P$-basic connection. Then there exists a twisted  basic Thom form $\tau$ on $E$.   Moreover,  the 
 fiber integration
$\pi_*: \Omega_{cv}(E, \mathcal{F}_E)[r] \rightarrow \Omega(M, \mathcal{O}_E, \mathcal{F})$
is a homotopy equivalence. A homotopy inverse of $\pi_*$ is the Thom map $\zeta_*: \Omega(M, \mathcal{O}_E, \mathcal{F}) \rightarrow \Omega_{ cv}(E,\mathcal{F}_E)[r]$
defined by $\zeta_*(\alpha) = \tau \wedge \pi^*\alpha$. \end{theorem}

\end{document}